\documentclass[12pt]{amsart}

\usepackage{latexsym,amssymb,amsthm,amsmath,amscd,multirow,graphics,multicol}
\usepackage{amssymb}

\theoremstyle{plain}  
\newtheorem{theorem}{Theorem}[section]

\newtheorem{proposition}{Proposition}[section]

\newtheorem{corollary}{Corollary}[section]

\newtheorem{example}{Example}[section]
\numberwithin{equation}{section}

\theoremstyle{remark}
\newtheorem{remark}{Remark}[section]

 \numberwithin{equation}{section}
\def\<{\left < }
\def\>{\right >}
\def\({\left ( }
\def\){\right )}
\def\e{\eqref}

\def\E4{\mathbb E^4 }

\def\d{\delta }
\def\with{with parallel mean curvature vector }

\def\x{\otimes}

\def\csch{\,{\rm csch\,}}
\def\sech{\hskip.012in{\rm sech\,}}
\def\e{\eqref}
\def\i{\hskip.01in {\rm i}\hskip.008in}

\begin{document}

\title[Submanifolds with parallel mean curvature vector]
{Submanifolds with parallel mean curvature vector in Riemannian and indefinite space forms}

\author[B. Y. Chen]{Bang-Yen Chen}

 \address{Department of Mathematics, 
	Michigan State University, East Lansing, Michigan 48824--1027, USA}
\email{bychen@math.msu.edu}

\begin{abstract}  A submanifold of a pseudo-Riemannian manifold is said to have {\it parallel mean curvature vector}  if the mean curvature vector field $H$ is parallel as a section of the normal bundle. Submanifolds with parallel mean curvature vector are important since they are critical points of some natural functionals.
 In this paper, we survey some classical and recent results on submanifolds with parallel mean curvature vector. Special attention is paid to the classification of space-like and Lorentz surfaces with parallel mean curvature vector in Riemannian and indefinite space forms. 
\end{abstract}


 \thanks{This is an invited survey article. The author thanks Prof. M. A. Al-Gwaiz, the editor-in-chief, for his invitation.}
\maketitle

\section{Introduction.}

A submanifold of a pseudo-Riemannian manifold is said to have {\it parallel mean curvature vector}  if the mean curvature vector field $H$ is parallel as a section of the normal bundle. A hypersurface of a Riemannian manifold has parallel mean curvature vector if and only if it has constant mean curvature. 
Trivially, every minimal submanifold of a pseudo-Riemannian manifold  has parallel mean curvature vector.  Furthermore, submanifolds with parallel second fundamental form have parallel mean curvature vector.

Submanifolds with parallel mean curvature vector are important since they are critical points of some natural functionals. In fact, a hypersurface of constant mean curvature in a Euclidean space is a solution to a variational problem. In particular, with respect to any volume-preserving variation of a domain $D$ in a Euclidean space the mean curvature of $M=\partial D$ is constant if and only if the volume of $M$ is critical. The condition of parallel mean curvature vector for submanifolds in higher dimensional Euclidean spaces are interest as well since, like its classical counterpart, are equivalent to a variational problem, namely, their Gauss map are harmonic map.

It was S. S. Chern who first suggested in the mid 1960s that the notion of parallel mean curvature vector as the natural extension of constant mean curvature for hypersurfaces.

In this paper, we survey some classical and recent results concerning submanifolds with parallel mean curvature vector in Riemannian manifolds as well as in pseudo-Riemannian manifolds.  In particular, we pay special attention to recent  progress on classification of space-like and Lorentz surfaces with parallel mean curvature vector in Riemannian and indefinite space forms. 

\section{Preliminaries.}

\subsection{Basic definitions, formulas and notations.}

Let $M$ be an $n$-dimensional  pseudo-Riemannian manifold isometrically immersed in a pseudo-Riemannian $m$-manifold  $\tilde M^m$. Let $\<\;\,,\;\>$ denote the inner product of $\tilde M^m$ as well as of $M$.
Denote by $\nabla$ and $\tilde \nabla$ the Levi-Civita connections on $M$ and $\tilde M^m$ respectively. Let  $h, D$ and $A$ be the second fundamental form, the normal connection,  and shape operator of $M$, respectively.

The Gauss and Weingarten formulas are given  by
\begin{align} & \label{2.1}\tilde \nabla_X Y=\nabla_X Y  + h(X, Y),\\& \label{2.2}
\tilde \nabla_X V  =-A_V X+D_X V\end{align}
for vector fields $X, Y$ tangent to $M$ and $V$ normal to $M$. For each $\xi\in T_x^{\perp}M$, the shape operator $A_{\xi}$ is a symmetric endomorphism of the tangent space $T_xM$ at $x\in M$. The shape operator and the second fundamental form are related by
\begin{align}\label{0.4} \<h(X,Y),\xi\>=\<A_{\xi}X,Y\>\end{align}
for $X,Y$ tangent to $M$ and $\xi$ normal to $M$.

Denote by $R$ and $\tilde{R}$ the Riemann curvature tensors of $M$ and $\tilde M^m$, respectively. Then the Gauss and Codazzi equations are given  by

\begin{align}&\label{2.3}\left< R(X,Y)Z,W\right>=\left< \right.\hskip-.02in \tilde
R(X,Y)Z,W\hskip-.02in \left.\right>+\left<h(X,W),
h(Y,Z)\right>\\& \hskip1.6in \notag -
\left<h(X,Z),h(Y,W)\right>,\\&\label{2.4} (\tilde R(X,Y)Z)^\perp=
(\bar\nabla_Xh)(Y,Z)-(\bar\nabla_Yh)(X,Z),
\end{align}
where $X,Y,Z,W$ are tangent vectors of $M$, $(\tilde R(X,Y)Z)^\perp$ is the normal component of
$\tilde R(X,Y)Z$,  and  $\bar \nabla h$ is defined by
\begin{equation}\label{2.5}(\bar \nabla_X h)(Y,Z) = D_X h(Y,Z) - h(\nabla_X Y,Z) - h(Y,\nabla_X Z).\end{equation}
When the ambient space $\tilde M^m$ is of constant curvature $c$, equations \e{2.3} and \e{2.4} of Gauss and Codazzi reduce to
\begin{align}&\label{2.6}\left< R(X,Y)Z,W\right>=c\{\<X,W\>\<Y,Z\>-\<X,Z\>\<Y,W\>\}\\& \hskip.6in \notag+\left<h(X,W), h(Y,Z)\right> - \left<h(X,Z),h(Y,W)\right>,\\&\label{2.7}
(\bar\nabla_Xh)(Y,Z)=(\bar\nabla_Yh)(X,Z),\end{align}
The mean curvature vector of $M$ in $\tilde M^m$ is defined by
\begin{align}&\label{2.8}H=\frac{{\rm trace}\, h}{\dim M}.\end{align}

The mean curvature vector $H$ is called parallel if we have $DH=0$ identically. A submanifold is said to have parallel second fundamental form if $\bar \nabla h=0$ holds identically. A submanifold with parallel second fundamental form is also known as a parallel submanifold.

A submanifold $M$ of a pseudo-Riemannian manifold is called {\it totally geodesic} if the second fundamental form vanishes identically.  It  is called {\it totally umbilical} if its second fundamental form satisfies $h(X,Y)=\<X,Y\>H$.
By an {\it extrinsic sphere} we mean a totally umbilical submanifold with nonzero parallel mean curvature vector. A circle is a 1-dimensional extrinsic sphere.

One-dimensional submanifolds with parallel mean curvature vector are nothing but
geodesics and circles. 

By a {\it CMC} surface of a pseudo-Riemannian 3-manifold, we mean a surface with nonzero constant mean curvature, i.e., a surface whose mean curvature vector $H$ satisfies $\<H,H\>=constant\ne 0$.

\subsection{Indefinite space forms.}

A pseudo-Riemannian manifold $\tilde M^m$ of constant sectional curvature is called {\it indefinite space form}. In particular, if $\tilde M^m$ is a Riemannian manifold of constant sectional curvature, it is called a Riemannian space form. An $m$-dimensional complete,  simply-connected Riemannian space form is isometric to an $m$-sphere $S^m$, a Euclidean $m$-space $\mathbb E^m$, or a hyperbolic $m$-space $H^m$ depending on the curvature is positive, zero, or negative. 

Indefinite space forms can be obtained as follows.
Let $\mathbb E_{t}^m$ denote the  pseudo-Euclidean $m$-space with the canonical pseudo-Euclidean metric of  index $t$ given by
\begin{align}\label{2.9}g_{0}= -\sum_{i=1}^{t} dx_{i}^{2} + \sum_{j=t+1}^{m}dx_{j}^{2}, \end{align}
where $(x_{1},\ldots,x_{m})$ is a rectangular coordinate system of  $\mathbb E_{t}^m$.   We put
\begin{align} &\label{2.10} S^k_s(x_0,c)=\big\{x\in \mathbb E^{k+1}_s| \<x-x_0,x-x_0\>={c}^{-1}>0\big\},\\&\label{2.11} H^k_s(x_0,-c)=\big\{x\in \mathbb E^{k+1}_{s+1}| \<x-x_0,x-x_0\>=-c^{-1}<0\big\},\end{align}
where $\<\;,\:\>$ is the indefinite inner product on $\mathbb E^{k+1}_t$. The $S^k_s(x_0,c)$ and $H^k_s(x_0,c)$ are complete pseudo-Riemannian manifolds with index $s$ of constant curvature $c$ and $-c$, respectively. The $S^k_s(x_0,c)$ and $H^k_s(x_0,c)$ are called {\it pseudo-Riemannian $m$-sphere} and the {\it pseudo-hyperbolic $m$-space}, respectively. 

We simply denote $S_s^k(x_0,c)$ and $H^k_s(x_0,-c)$ by $S^k_s(c)$ and $H^k_1(-c)$, respectively, when $x_0$ is the origin.
The pseudo-Riemannian manifolds $\mathbb E^n_s, S^m_s(c)$ and $H^m_s(-c)$ are  {\it indefinite space forms}. 
In general relativity, the Lorentzian manifolds $\mathbb E^n_1,S^n_1(c)$ and $ H^n_1(-c)$ are known as the  Minkowski, de Sitter, and anti-de Sitter spaces (or space-times), respectively. These Lorentzian manifolds have constant sectional curvatures which are called {\it Lorentzian space forms}.

A  vector $v$ is called {\it space-like} (resp., {\it time-like}) if  $\<v,v\>>0$ (resp., $\<v,v\><0$). A vector $v$ is called  {\it light-like} if it is  nonzero and it satisfies $\<v,v\>=0$. A submanifold of pseudo-Riemannian manifold is called space-like if each nonzero tangent vector is space-like.

 A curve $z=z(t)$  in $R^m_s(c)$ defined on an open interval $I$ is called a {\it null curve} (resp., a {\it space-like curve} or a {\it time-like curve}) if its velocity vector $z'(t)$ is a light-like (resp., space-like or time-like) at each $t\in I$. Two curves $z,w$ defined on $I$ are called orthogonal if $\<z'(t),w'(t)\>=0$ at each $t\in I$.  
 
The {\it light cone} $\mathcal{LC}(c_0)$ with vertex $c_0$ in $\mathbb E^{m}_s \, (1\leq s<m)$ is defined by
\begin{align}\notag\mathcal{LC}(c_0)=\big\{x \in
\mathbb E^{m}_s\,:\, \<x-c_0,x-c_0\>=0\big\}.\end{align}
When $c_0$ is the origin, $\mathcal{LC}(c_0)$ is simply denoted by $\mathcal{LC}$. 

\section{Gauss map and parallel mean curvature vector.}

Harmonic maps $\psi:(M,g)\to (N,h)$ between Riemannian manifolds are critical points of the energy functional:
$$E(\psi)=\frac{1}{2}\int_M |d\psi|^2v_g,$$
where  $v_g$ is the volume element of $(M,g)$. The Euler-Lagrange equation associated with $E$ is given by vanishing of the tensor field $\tau(\psi)={\rm trace}(\nabla d\psi)$.
If the map $\psi$ is isometric, then up to a nonzero constant the tension field reduces to the mean curvature vector field (cf. \cite[1964]{ES}).

For an isometric immersion $\phi: M\to \mathbb E^m$  of  an oriented  Riemannian $n$-manifold into a Euclidean m-space, the Gauss map
$$G : M\to G^R(m-n,m)$$ of $\phi$ is a smooth map which carries a point $x\in M$ into the oriented $(m-n)$-plane in $\mathbb E^m$, which is
obtained from the parallel translation of the normal space of $M$ at $x$ in $\mathbb E^m$,  
where $G^R(m-n,m)$ denotes the  Grassmannian manifold consisting of oriented 
$(m-n)$-planes in $\mathbb E^m$.

Ruh and Vilms \cite[1970]{RV}  characterized submanifolds of Euclidean spaces with
parallel mean curvature vector  as follows.

\begin{theorem} A submanifold $M$ of a Euclidean $m$-space
$\mathbb E^m$ has parallel mean curvature vector if and only if its Gauss map $G$ is harmonic.
\end{theorem}

\section{Construction and existence of submanifolds with parallel mean curvature vector.}

Reckziegel \cite[1974]{Re} proved that if $M$ is a compact submanifold of a
manifold $N$ such that the restriction of the tangent bundle $TN$ to $M$, $TN|_M$, has a
metric  $g$ and $\eta$ is a nonzero normal vector field of constant length in $(TM)^\perp$, then $g$ can be extended to a Riemannian metric on $N$ such that $M$ is an extrinsic sphere with parallel mean curvature vector $\eta$. This result shows  that extrinsic spheres have no special topological properties.

In order to construct many more non-trivial examples of totally umbilical submanifolds with parallel mean curvature vector in Riemannian manifolds, the author introduced in \cite[page 66]{c1981} the notion of twisted products as follows:

Let $B$ and $F$ be Riemannian manifolds with Riemannian metrics $g_B$ and $g_F$, respectively, and $f$  a positive differentiable function on $B\times F$. Consider the product manifold
$B\times F$ with its projection $\pi_B:B\times F\to B$ and $\pi_F:B\times F\to F$. The {\it  twisted product} $B\times_f F$ is the manifold $B\times F$ equipped with the Riemannian
structure such that
\begin{align} \notag ||X||^2=||\pi_B{}_*(X)||^2+f^2 ||\pi_{F}{}_*(X)||^2\end{align} for any  vector $X$ tangent to $B\times_f F$. Thus, we have $g=g_B+f^2 g_F$.

In \cite[1981]{c1981}, the author proved that every
Riemannian manifold can be embedded in some twisted product Riemannian manifold as an extrinsic sphere. More precisely, he proved the following result.

\begin{theorem} Let $M=B\times_f F$ be a twisted product of two Riemannian manifolds $B$ and $F$. Then  we have:
\vskip.04in

{\rm (1)} For each $b\in B$, the fiber $F_b=\{b\}\times F$ is a totally umbilical submanifold in $M$ with $f^{-1}\nabla f$ as its mean curvature vector, where $\nabla f$ is the gradient of $f$.
\vskip.04in

{\rm (2)} Fibers have constant mean curvature if and only if  $||\nabla (\ln f)||$ is a function of $B$.
\vskip.04in

{\rm (3)} Fibers have parallel mean curvature vector if and only if $f$ is the product of two positive functions $\lambda(b)$ and $\mu(p)$ of $B$ and $F$, respectively.
\end{theorem}

Using the method of equivariant differential geometry, W. T. Hsiang, W. Y. Hsiang and Sterling \cite[1985]{HHS}  obtained the following.

\begin{theorem} We have:
\vskip.04in

{\rm (a)} There exist infinitely many codimension two embeddings of distinct knot types of $S^{4k+1}$ into $S^{4k+3}(1)$ with parallel mean curvature vector of arbitrarily small constant length. 
\vskip.04in

{\rm (b)} There exist infinitely many codimension two embeddings of distinct knot types of the Kervaire exotic sphere  $\Sigma^{4k+1}_0$ into $S^{4k+3}(1)$ with parallel mean curvature vector having
length of arbitrarily small constant value. 
\vskip.04in

{\rm (c)} There exist infinitely many constant mean curvature embeddings of $(4k-1)$-dimensional generalized lens spaces into $S^{4k+1}(1)$. 
\end{theorem}

\section{Surface with parallel mean curvature vector in Riemannian space forms.}

Ferus \cite[1971]{Ferus1} and Ruh \cite[1971]{Ruh} determined  closed surfaces of genus zero with parallel mean curvature vector in a Euclidean space as follows.

\begin{theorem} Let $M$ be a closed oriented  surface of genus zero in $\mathbb E^m$. If $M$ has parallel mean curvature vector, then $M$ is contained in a hypersphere of $\mathbb E^m$ as a minimal surface.\end{theorem}

For surfaces of constant Gaussian curvature, Chen and  Ludden \cite[1972]{CL} and  Hoffman  \cite[1973]{Ho} proved the following result.

\begin{theorem}\label{CL} Minimal surfaces of a small hypersphere, open pieces of the product of two plane circles, and open pieces of a circular cylinder are the only non-minimal surfaces in a Euclidean space with parallel mean curvature vector and with constant Gaussian curvature.
\end{theorem}

Surfaces in $\mathbb E^4$ with parallel mean curvature vector were determined by Hoffman in his doctoral thesis [Stanford University, 1971]. 

The complete classification of  surfaces in Euclidean $m$-space, $m\geq 4$, with parallel mean
curvature vector was obtained by Chen \cite[1973]{c} and Yau \cite[1974]{Yau}.

\begin{theorem}  A surface $M$ of a Euclidean $m$-space $\mathbb E^m$ has parallel mean curvature vector if and only if  it is one of the following surfaces:
\vskip.04in

{\rm (a)}  a minimal surface of $\mathbb E^m$; 
\vskip.04in

{\rm (b)}  a minimal surface of a hypersphere of $\mathbb E^m$;
\vskip.04in

{\rm (c)}  a surface of $\mathbb E^3$ with constant mean curvature;
\vskip.04in

{\rm (d)}  a surface of constant mean curvature lying in a hypersphere of an affine 4-subspace of $\mathbb E^m$. \end{theorem}

Similar results hold for surfaces with parallel mean curvature vector in spheres and in real hyperbolic spaces as well (cf. \cite[1973]{c1973}).

For  a compact surface $M$ with positive constant Gaussian curvature,  Enomoto \cite[1985]{En}  proved that if  $\phi :M\to \mathbb E^m$ is an isometric embedding with constant  mean curvature and flat normal  connection, then $\phi(M)$ is a round sphere in an affine 3-subspace of $\mathbb E^m$.
 Enomoto also proved that if $f:M\to \mathbb E^{n+2}$ is an isometric embedding of a compact
Riemannian $n$-manifold, $ n\geq 4$, with positive constant sectional curvature and with constant mean curvature, then $\phi(M)$ is a round $n$-sphere in a hyperplane of $\mathbb E^{n+2}$.

It remains as an open problem to  completely classify submanifolds of dimension $\geq 3$ with parallel mean curvature vector in Riemannian space forms.

\section{Surfaces with parallel normalized mean curvature vector.} 

In  \cite[1980]{c1980}, the author  extended the notion of submanifolds \with  to the notion of  submanifolds with parallel normalized mean curvature vector. He called a submanifold $M$ to have {\it parallel normalized mean curvature vector} if locally there exists a unit parallel vector field $\xi$ on $M$ which is parallel to the mean curvature vector $H$, that is, $H=\alpha \xi$ for
some unit parallel normal vector field $\xi$.

Obviously, every submanifold with nonzero parallel mean curvature vector   has parallel normalized mean curvature vector.  However, the condition to have parallel normalized mean curvature vector  is much weaker than the condition to have parallel mean curvature vector. For instance, every hypersurface in a Riemannian manifold always has parallel normalized mean curvature vector.

We have the following result from \cite[1980]{c1980} for surfaces with parallel normalized mean curvature vector.

\begin{theorem} \label{T:6.1} Let $M$ be an analytic surface in a complete, simply-connected Riemannian space form $R^m(c)$ of constant curvature $c$. If $M$ has parallel normalized mean curvature vector, then either $M$ lies in a hypersphere of $R^m(c)$ as a minimal surface or $M$ lies in a 4-dimensional totally geodesic submanifold of $R^m(c)$.
\end{theorem}

As applications of Theorem \ref{T:6.1} he also obtained in \cite[1980]{c1980} the following.

\begin{theorem} \label{T:6.2} Let $M$ be a flat analytic surface in a Euclidean $m$-space. If $M$ has parallel normalized mean curvature vector, then $M$ is one of the following:

\vskip.04in
{\rm (1)} a flat minimal surface of a hypersphere of $\mathbb E^m$,
\vskip.04in

{\rm (2)} an open piece of the product of two plane circles, or
\vskip.04in

{\rm (3)} a developable surface in a linear 3-subspace of $\mathbb E^m$.
\end{theorem}

\begin{theorem} \label{T:6.3} Let $M$ be a  Riemann sphere in a Euclidean $m$-space $\mathbb E^m$. If $M$ has parallel normalized mean curvature vector, then either $M$ lies in a hypersphere of $\mathbb E^m$ as a minimal surface or $M$ lies in  a linear 3-subspace of $\mathbb E^m$.
\end{theorem}

\begin{theorem} \label{T:6.4} Let $M$ be a closed analytic surface in a Euclidean $m$-space $\mathbb E^m$. If $M$ has constant Gauss curvature and parallel normalized mean curvature vector, then either $M$ is a minimal surface of a hypersphere of $\mathbb E^m$ or $M$ is the product of two plane circles.
\end{theorem}

Every surface in a Euclidean 3-space has parallel normalized mean curvature vector. Moreover, there exist abundant examples of surfaces which lie fully in a Euclidean 4-space with parallel normalized mean curvature vector, but  not with
parallel mean curvature vector.

\section{Finite type submanifolds and parallel mean curvature vector.}

A smooth map $\phi$ from a (pseudo) Riemannian manifold $M$ into a (pseudo) Euclidean space $\mathbb E^m$ is said to be {\it of finite type} if $\phi $ is decomposable as a finite sum of $\mathbb E^m$-valued eigenfunctions of the Laplacian $\Delta$ on $M.$ In particular, a map $\phi $ is of 2-type if and only if $$\phi = \phi_0 + \phi_1+\phi_2 ,$$ where $\phi_0 =$ constant, $\, \phi_1,\phi_2 \neq$ constant and $\Delta \phi_i = \lambda_i \phi_i$ with $\lambda_1\ne \lambda_2$. When one of $\lambda_1,\lambda_2$ is zero, then $\phi$ is said to be of {\it null 2-type}.

If  $\phi \colon M\to \mathbb E^m$ is an isometric immersion and has a decomposition into a finite sum of eigenfunctions, then $M$ is said to be of finite type in $\mathbb E^m$ (via $\phi$) (see \cite[1984]{c1984} and \cite[1996]{report} for details).

In \cite[1988]{CLue}, the author and Lue investigated the relationship  between 2-type submanifolds and submanifolds with parallel mean curvature vector. Consequently, they obtained the following.

\begin{theorem} Let $M$ be a 2-type submanifold of a Euclidean $m$-space $\mathbb E^m$. If $M$ has parallel mean curvature vector, then either $M$ is spherical $($i.e., $M$ lies in a hypersphere of $\mathbb E^m)$ or $M$ is of null 2-type.
\end{theorem}

\begin{theorem} Let $M$ be a closed 2-type surface of a Euclidean $m$-space $\mathbb E^m$. Then $M$ has parallel mean curvature vector if and only if $M$ is the product of two plane circles with different radii.
\end{theorem}

\begin{theorem} Let $M$ be a  surface of a Euclidean $m$-space $\mathbb E^m$ with parallel mean curvature vector. Then $M$ is of 2-type if and only if $M$ is an open portion of the product surface of two plane circles with different radii or $M$ is an open portion of a circular cylinder.
\end{theorem}

In \cite[1994]{Li}, Li proved that a surface in $\mathbb E^m$ with parallel normalized mean curvature vector is of null 2-type if and only if it is an open portion of a circular cylinder. When $m=4$ this result was due to the author \cite[1989]{c1989}.

Dursun classified in \cite[2005]{Dur}  3-dimensional null 2-type submanifolds in ${\bold E^5}$ with parallel normalized mean curvature vector under the conditions that the shape operator in the direction of the mean curvature vector has at most two distinct eigenvalues and that the squared norm of the second fundamental form is constant.

Null 2-type surfaces in the Minkowski 4-space $\mathbb E^4_1$ with parallel mean curvature vector were studied by Kim in \cite[1997]{YHKim}. He showed that  if $M$ is a null 2-type surface in $\mathbb E^4_1$ with unit parallel mean curvature vector. Then $M$ is an open portion of a $B$-scroll in $\mathbb E^3_1$, $\mathbb E^1_1\times S^1(r)$, $S_1^1(r)\times {\bf R}$, $H^1(r)\times {\bf R}$ or an extended $B$-scroll in $\mathbb E^4_1$. As a consequence, the only null 2-type surface  with unit parallel mean curvature vector lying fully in $\mathbb E^4_1$ is an open part of an extended $B$-scroll.

\section{Tensor product immersions with parallel mean curvature vector.}

Let $V$ and $W$ be two vector spaces over the field of real  numbers $\bf R$. Then we have the notion of the tensor product $V\otimes W$ defined as follows (cf. \cite[1993]{c1993}). 

If $V$ and $W$ are inner product spaces with their inner products given respectively by $\<\; ,\; \>_V$ and $\<\; ,\; \>_W$, then  $V\otimes W$ is also an inner product space with  inner product defined  by $$ \<v\otimes w,z\otimes
y\>=\<v,z\>_V\<w,y\>_W \quad  \forall v\otimes w, z\otimes y \in V\otimes W.$$  
Let $\mathbb E^m$ denote the  Euclidean $m$-space with the canonical Euclidean inner product. Then, with respect to
the inner product defined above,   $\mathbb E^m \otimes \mathbb E^{m'}$ is isometric to $\mathbb E^{mm'}.$   By applying
this  algebraic notion, we have the notion of  {\it tensor product map\/}
$$f_1\otimes f_2 : M \rightarrow \mathbb E^m \otimes \mathbb E^{m'}\equiv \mathbb E^{mm'}$$
 associated with any two  maps $f_1:M\to \mathbb E^m$ and $f_2:M\to \mathbb E^{m'}$ of a given Riemannian manifold $(M,g)$ defined as follows:
$$(f_1\otimes f_2)(x)=f_1(x)\otimes f_2(x) \in  \mathbb E^m \x \mathbb E^{m'},\quad  \forall x\in M.$$
An immersion $f:M\rightarrow \mathbb E^m$ is said to be {\it spherical\/} if $f(M)$ is contained in a hypersphere of $\mathbb E^m$ centered at the origin of $\mathbb E^m$.

Denote by ${\mathcal S}(M)$ the set of all spherical immersions from a  Riemannian manifold $(M,g)$ into Euclidean spaces. Then  $\otimes$ is a  binary operation on ${\mathcal S}(M)$.   Hence, if   $f_1:M\to \mathbb E^m$ and $f_2:M\to \mathbb E^{m'}$ are immersions belonging to  ${\mathcal S}(M)$,  their tensor product map $f_1\otimes f_2:M\to \mathbb E^m\otimes \mathbb E^{m'}$ is an immersion in ${\mathcal S}(M)$,  called the {\it tensor product immersion\/} of $f_1$ and $f_2$. 
 
 Similarly, if $f_i:M\rightarrow \mathbb E^{m_i},\; i=1,\ldots,t,$ are $t$ immersions belonging to ${\mathcal S}(M)$,  one has the tensor product immersion $f_1\x\cdots\x f_t$ of $f_1,\ldots,f_t$. 
 
In \cite[1991]{c1991}, the author studied tensor product immersions which have parallel mean curvature vector. He obtained the following.

\begin{theorem} Let $k \ge 2$ and $f_i\colon M \rightarrow S^{m_i -1}(r_i)$, $i=1,\cdots,k$, be $k$ spherical isometric immersions from a Riemannian manifold $(M,g)$. Then the tensor product immersion $f_1 \otimes \cdots \otimes f_k$ has parallel mean curvature vector if and only if $k = 2$, $r_1=r_2$ and $f_1,f_2$ are totally geodesic.
\end{theorem}

\section{Graph with parallel mean curvature vector.}

Let $f:M\to N$ be a smooth map, where $M$ and $N$ are Riemannian manifolds of dimensions $m$ and $n$ with Riemannian metrics $g$ and $h$, respectively. The graph of $f$, $\Gamma_f=\{(x,f(x)): x\in M\}$, is a submanifold of $M\times N$. Take the product metric on $M\times N$, and the induced metric on $\Gamma_f$. 

Salavessa \cite[1989]{Sa} proved the following theorem. 

\begin{theorem} Let $f:M\to N$ be a smooth map.
If the graph $\Gamma_f : M\to M\times N$ is an immersion with parallel mean curvature vector $H_{\Gamma_f}$, then, for each oriented compact domain $D\subset M$, we have the isoperimetric inequality $$c\leq \frac{1}{m}\frac{A(\partial D)}{V(D)},$$ where $c=||H_{\Gamma_f}||_{g\times h}$ $(c$ a constant$\,)$ and where $V(D)$ $($resp. $A(\partial D))$ is the volume of $D$ $($resp. the area of $\partial D)$ relative to the metric $g$. 

In other words, if $(M,g)$ is an oriented Riemannian manifold, then $$|| H_{\Gamma_f}|| _{g\times  h}\leq\frac{1}{m}h(M),$$ where $h(M)$ is the Cheeger constant of $M$. In particular, if $(M,g)$ has zero Cheeger constant, then $\Gamma_f$ is a minimal submanifold of $M\times N$. \end{theorem}

\section{Homogeneous submanifolds with parallel mean curvature vector.}

 Olmos investigated in \cite[1994]{O1} and \cite[1995]{O2} homogeneous
submanifolds with parallel mean curvature vector in a Euclidean space. He obtained the following.

\begin{theorem} We have:
\vskip.04in

{\rm (a)} If $M$ is a compact homogeneous submanifold of a Euclidean space with parallel mean curvature vector which is not minimal in a sphere, then $M$ is an orbit of the isotropy representation of a simple symmetric space; 
\vskip.04in

{\rm  (b)} A homogeneous irreducible submanifold of Euclidean space with parallel mean  curvature vector is either minimal, or minimal in a sphere, or an orbit of the isotropy  representation of a simple Riemannian symmetric space.
\end{theorem}

\section{Slant surfaces with parallel mean curvature vector.}

Let $M$ be a Riemannian $n$-manifold isometrically immersed in an almost Hermitian manifold
$\tilde M$ with (almost) complex structure $J$ and almost Hermitian metric $g$. 
For any vector $X$ tangent to $M$ we put
$$JX=PX+FX,$$
where $PX$ and $ FX$ are the tangential and the normal components of $JX$, respectively.  Thus, $P$ is an endomorphim of the tangent bundle $TM$ and $F$ a normal-bundle-valued 1-form on $TM$.

For any nonzero vector $X$ tangent to $M$ at a point $x \in M$, the angle $\theta (X)$  between $JX$ and the tangent space $T_{x}M$ is called the {\it Wirtinger angle\/} of $X$. In the following  we call an immersion $\phi :M \rightarrow \tilde M$ a {\it slant immersion\/} if the Wirtinger angle $\theta (X)$ is constant (which is independent of the choice of $x \in M$ and of $X \in T_{x}M$). Complex and totally real immersions are slant immersions with $\theta = 0$ and $\theta = \pi /2$,
respectively. Slant submanifolds of an almost Hermitian manifolds are characterized by the condition: $\, P^{2}=\lambda I,\,$ for some real number $\lambda \in [-1,0],$ where $I$ denotes the identity transformation of the tangent bundle $TM$ of the submanifold $M$.
 The Wirtinger angle of a slant immersion is called the {\it slant angle}. A slant submanifold is said to be {\it proper\/} if it is neither complex nor totally real  (see \cite[1990]{c1990} for details). 

Let ${\bf C}^m,\,CP^m(4\rho)$ and $CH^m(-4\rho)$  denote the flat complex Euclidean $m$-plane, the complex projective $m$-space of constant holomorphic sectional curvature $4\rho$ and the complex hyperbolic $m$-space with constant holomorphic sectional curvature $-4\rho$, respectively.

In \cite[1990, page 50]{c1990}, the author classified slant surfaces in ${\bf C}^2$ with parallel mean curvature vector in the following.

\begin{theorem}\label{T:7,1} Let $M$ be a surface in ${\bf C}^2$. Then $M$ is a slant surface with parallel mean
curvature vector if and only if $M$ is one of the following surfaces:

{\rm (a)} an open portion of the product surface of two plane circles, or

{\rm (b)} an open portion of a circular cylinder which is contained in a hyperplane of ${\bf C}^2$, or

{\rm (c)} a minimal slant surface in  ${\bf C}^2$.
\end{theorem}

Also,  the  author proved  in \cite[1998]{c1998} the following result.

\begin{theorem}  The squared mean curvature $H^2=\<H,H\>$ and  the Gauss curvature $K$ of  a proper slant surface $M$ in a complex space form $\tilde M^2(4\epsilon)$ of constant holomorphic sectional curvature $4\epsilon$ satisfy the following basic inequality:
\begin{align}\label{7.1} H^2\geq 2K-2(1+3\cos^2\theta)\epsilon,\end{align}
\noindent  where $\theta$ is the slant angle.  

Moreover,  the equality sign of \eqref{7.1} holds identically if and
only if, with respect to some suitable adapted orthonormal frame $\{e_1,e_2,e_3,e_4\}$, the shape operator of $M$  takes the following special form:
\begin{align}\label{7.2}A_{e_3}=\begin{pmatrix} 3\lambda & 0\\ 0 &
\lambda\end{pmatrix},\quad A_{e_4}=\begin{pmatrix} 0 &
\lambda\\ \lambda & 0\end{pmatrix}\end{align}
for some function $\lambda$. \end{theorem}

The author proved  in \cite[1998]{c1998} that there do not exist  proper slant surfaces satisfying the equality case of inequality \eqref{7.1} for $\epsilon>0$. A proper slant surface in a flat K\"ahler surface satisfies the equality of inequality \eqref{7.1} if and only if it is totally geodesic. Furthermore, he proved that a proper slant surface in the complex
hyperbolic plane $CH^2(-4)$ satisfying the equality case of  inequality \eqref{7.1} identically is a surface of constant Gaussian curvature $-{2\over 3}$ with slant angle $\theta=\cos^{-1}({1\over 3})$.
The immersion of this slant surface in $CH^2(-4)$ is rigid and has parallel mean curvature vector.

A submanifold $M$ of a pseudo-Riemannian Sasakian manifold $(\tilde M,g,\phi,\xi)$ is called {\it contact $\theta$-slant} if the structure vector field $\xi$ of $\tilde M$ is tangent to $M$ at each point; moreover, for each unit vector $X$ tangent to $M$ and orthogonal to the characteristic vector field $\xi$ at $p\in M$, the angle $\theta(X)$ between $\phi(X)$ and $T_pM$ is independent of the choice of $X$ and $p$.

Let $H^{2m+1}_1(-1)\subset {\bf C}^{m+1}_1$ denote the anti-de Sitter space-time and $\pi\colon\;
H^{2m+1}_1(-1)$ $\to  CH^m(-4)$ the corresponding totally geodesic fibration. Then each $n$-dimensional proper $\theta$-slant submanifold $M$ in $CH^{m}(-4)$ lifts to an $(n+1)$-dimensional
proper contact $\theta$-slant submanifold $\pi^{-1}(M)$ in $H^{2m+1}_1(-1)$ via $\pi$.

Conversely, every proper contact $\theta$-slant submanifold of $H^{2m+1}_1(-1)$ projects to a proper $\theta$-slant submanifold of $CH^{m}(-4)$ via $\pi$ (see \cite[2000]{CT} for details).

The contact slant representation of the unique proper slant surface in $CH^2(-4)$ which satisfies the equality case of \eqref{7.1}  in $H^5_1(-1)\subset C^3_1$ has been determined by the author and Tazawa in \cite[2000]{CT}. In fact, they proved that, up to rigid motions of ${\bf C}^3_1$, the contact slant representation of this surface of constant curvature $-\frac{2}{3}$ and with slant angle $\theta=\cos^{-1}({1\over 3})$ is given by
\begin{equation}\begin{aligned} \label{7.3} &z(u,v,t)=e^{it}\Big(1+\frac32\Big(\cosh
\sqrt{2\over3}v-1\Big)+\frac{u^2}6e^{-\sqrt{2\over3}
v}-i{u\over\sqrt{6}}
(1+e^{-\sqrt{2\over3}v}),\\&\hskip.2in  \frac u3\Big(
{1+2e^{-\sqrt{2\over3}v}}\Big)+\frac
i{6\sqrt{6}}e^{-\sqrt{2\over3}v}\Big(\Big(e^{\sqrt{2\over3}v}
-1\Big)\Big(9e^{\sqrt{2\over3}v}-3\Big)
+2u^2\Big),\\ &\hskip.3in \frac u{3\sqrt{2}}
\Big(1-e^{-\sqrt{2\over3}v}\Big)+\frac
i{12\sqrt{3}} \Big(6-15 e^{-\sqrt{2\over3}v}+
9e^{\sqrt{2\over3}v}
+2e^{-\sqrt{2\over3}v}u^2\Big)\Big).
\end{aligned}\end{equation}
This slant surface with parallel mean curvature vector and with constant Gauss curvature in $CH^2(-4)$ is named as Chen's surface by Kenmotsu and Zhou  in \cite[2000]{KZ}.

\section{Surfaces \with in complex space forms.}

In the form of an overdetermined system of differential equations, Ogata  studied in \cite[1995]{Og} non-minimal immersions of a surface $M$ into the complex projective plane $CP^2(4)$ with parallel curvature vector, showing a method for the local construction of such immersions. 

In \cite[2000]{KZ}, Kenmotsu and Zhou investigated surfaces with nonzero parallel mean curvature vector in 2-dimensional complex space forms  by solving Ogata's overdetermined system without assuming the slant condition. More precisely, they claimed the following results.

\begin{proposition}\label{T:7.2}  Let $\phi:M\to CP^2(4)$ be an isometric immersion of a surface $M$ into $CP^2(4)$ with nonzero parallel mean curvature vector. Then $M$ is flat and $\phi$  is totally real.\end{proposition}

\begin{proposition} There exists a one-parameter family of isometric surfaces in ${\bf C}^2$ that preserves the length of the mean curvature vector and the Wintinger angle. Each
surface of the family has the following properties:
\vskip.04in

{\rm (1)} nonzero parallel mean curvature vector,
\vskip.04in

{\rm (2)} nonconstant Wirtinger angle,
\vskip.04in

{\rm (3)} the Riemannian metric isometric to a rotational surface in ${\mathbb E}^3$.
\vskip.04in

Conversely, let  $\phi:M\to {\bf C}^2$ be an isometric immersion with parallel mean curvature
vector with  $|H| = 2b > 0$. Then, either $\phi$  is totally real or locally congruent
to an element of the family above. \end{proposition}

\begin{proposition} \label{T:7.4} There exists a one-parameter family of surfaces in
$CH^2(-12b^2)$ that preserves the length of the mean curvature vector and the Wirtinger
angle. Each surface of the family has the following properties:
\vskip.04in

{\rm (1)} parallel mean curvature vector whose length is equal to 2b,

\vskip.04in
{\rm (2)} a proper slant surface,
\vskip.04in

{\rm (3)} the Riemannian metric isometric to a rotational surface in $\mathbb E^3$.

\vskip.04in

Conversely, let  $\phi:M\to CH^2(4\rho)$  be an isometric immersion with parallel
mean curvature vector with $|H|=2b>0$. Then we have $\rho=-3b^2$  and one of the
following holds:

\vskip.04in
{\rm (1)} $\phi$  is totally real,
\vskip.04in

{\rm (2)} $\phi$  is locally congruent to Chen's surface  if $M$ is of constant Gaussian 
curvature,
\vskip.04in

{\rm (3)} $\phi$  is locally congruent to a surface of the family above  if $M$  is not of constant
Gaussian curvature.
\end{proposition}

However, as pointed out by Hirakawa in \cite[2006]{Hi},  there exists a gap in the argument used by Ogata in order to derive his overdetermined system in \cite[1995]{Og}. In fact, Hirakawa observed that some special coordinates that are used in the proof of Ogata's paper are not available in general. As a consequence, the results obtained by Kenmotsu and Zhou in \cite[2000]{KZ} given above only hold for the restricted class of surfaces for which such special coordinates exist. 

In \cite[2006]{Hi} by Hirakawa partially fills the gap by classifying all surfaces with nonzero parallel mean curvature vector in a two-dimensional complex space form under the additional assumption that they have constant Gaussian curvature. In particular, he obtained explicit new examples of such surfaces. More precisely, he obtained the following.

\begin{theorem} Let $M$ be an oriented two-dimensional Riemannian manifold, and $\phi$ be an isometric immersion of $M$ into a two-dimensional complex space form with nonzero parallel mean curvature vector of length $2b$. If $M$ has constant Gaussian curvature $K$, then $K$ is equal to $-2b^2,0$ or $4b^2$.

\vskip.04in
{\rm (1)} When $K= -2b^2$, $\phi(M)$ is in $CH^2(-12b^2)$ and it is part of 

\vskip.04in
\noindent {\rm (1.1)} the slant surface of Chen defined by \e{7.3}, or
\vskip.04in

\noindent {\rm (1.2)} one of the surfaces described in Example 7.1 below. 
\vskip.04in

{\rm (2)} When $K=0$, we have
\vskip.04in

\noindent  {\rm (2.1)} $\phi$ is totally real and $\phi(M)$  is part of a product of two circles in $CP^2(4\rho)$, where $\rho>0$, or
\vskip.04in

\noindent {\rm (2.2)} $\phi(M)$  is part of a cylinder or a product of two circles in ${\bf C}^2$, or 

\vskip.04in
\noindent {\rm (2.3)} $\phi$ is totally real and $\phi(M)$ is part of a plane, a cylinder or a product of two circles in $CH^2(-4\rho)$, where $0<\rho\leq 2b^2$. 
\vskip.04in

{\rm (3)} When $K=4b^2$, $\phi(M)$  is part of a round sphere in a hyperplane in ${\bf C}^2$.
\end{theorem}

\begin{example} {\rm Let $b$ be a positive real number. Let us consider $CH^2(-12b^2)$ as a quotient space $S/S^1$, where 
\begin{align}\notag &S=\big\{(u_0,u_1,u_2)\in {\bf C}^3: -|u_0|^2+|u_1|^2+|u_2|^2=-1\big\},\;  \\&\notag S^1=\{e^{it}:t\in {\bf R}\}.\end{align}

 Denote by $\pi:S\to CH^2(-12b^2)$ the projection. Let $M$ be the 2-dimensional Riemannian manifold $(0,\frac{\pi}{2})\times {\bf R}$ in the real 2-plane ${\bf R}^2$ equipped with the metric
$$ds^2_M=\frac{2}{b^2\sin^22u}(du^2+dv^2),$$
where $(u,v)$ is the standard coordinate of ${\bf R}^2$. Then $M$ has constant Gaussian curvature $-2b^2$ and is complete. 
 
 Define two maps $X_\pm: M\to S$  by
$${}^tX_\pm (u,v)=\frac{\csc 2u}{2 \sqrt{1+\cos^2u}}\begin{pmatrix} 1&1&1\\ 0&1&1 \\1&0&1\end{pmatrix}\begin{pmatrix} e^v(3\pm 2\sqrt{2}\cos^2 u+e^{-2iu})\\ e^{-v}(3\mp 2\sqrt{2}\cos^2u+2^{2iu})\\ \pm 4\cos^3u+2\sqrt{2}i\sin u\end{pmatrix},
$$
  and define $x_{\pm}$  by $\pi\circ X_\pm$. Then $x_\pm:M\to CH^2(-12b^2)$ are isometric embeddings
with parallel mean curvature vector of length $2b$.
}\end{example}

Hirakawa also obtained the following.

\begin{corollary} Let $M$ be an oriented two-dimensional Riemannian manifold, and $\phi$ an isometric immersion of $M$ into a two-dimensional complex space form with nonzero parallel mean curvature vector of length $2b$. If $M$ is homeomorphic to a sphere, then $\phi(M)$ is a round sphere in a hyperplane in ${\bf C}^2$.
\end{corollary}

\section{Periodicity of plane in complex space forms.}

Let $\phi:{\mathbb E}^2\to \tilde M^2(4\rho)$ be an isometric immersion of a Euclidean plane ${\mathbb E}^2$ into a  complex space form $\tilde M^2(4\epsilon)$.
Assume that the immersion $\phi$ has nonzero parallel mean curvature vector.

When the ambient space $\tilde M^2(4\rho)$ is the complex Euclidean plane ${\bf C}^2$, Theorem \ref{CL} implies that $\phi$ must be either doubly periodic or singly periodic. 

When  the ambient space $\tilde M^2(4\rho)$ is $CP^2(4\rho)$,  the main result of \cite[1995]{DT} implies that $\phi$ is doubly periodic.

In \cite[2004]{Hir}, Hirakawa studied the periodicity of plane with nonzero parallel mean curvature vector in $CH^2(-4\rho)$. Now, we explain the result of Hirakawa as follows.
Kenmotsu and Zhou proved in \cite[2000]{KZ} that an isometric immersion $\phi:{\mathbb E}^2\to CH^2(-4\rho)$ is uniquely determined by two real numbers $b$ and $t$ which are defined by the second fundamental form and satisfy $0<b$ and $0\leq t<2\pi$, up to holomorphic isometries of $CH^2$. The real number $b$ is  half of the length of the mean curvature vector.  The other constant $t$ is determined from $\phi$ as follows: 

For the given $\phi:{\mathbb E}^2\to CH^2(-4\rho)$,  let $e_4$  be the unit normal vector field  orthogonal to $e_3=H/|H|$  and compatible to the orientation of $CH^2(-4\rho)$. Let $\{e_1,e_2\}$ be an oriented orthonormal frame field on $\mathbb E^2$. Denote
by $h^r_{ij}\,(r=3,4; i,j=1,2)$ the coefficients of the second fundamental form  with respect to $e_1,e_2,e_3,e_4$, 

The two quadratic differentials 
$$\varphi_r=(h^r_{11}-h^r_{22}-\sqrt{-1}\, h^r_{12})dz^2,\; r=3,4,$$
are globally defined on $\mathbb E^2$. The ratio $\varphi=\varphi_3/\varphi_4$ is constant and 
$$\Psi=\frac{\sqrt{1+\rho/2b^2} (\sqrt{-1}+\varphi)}{\sqrt{-1}-\varphi}$$ has absolute value 1. The real number $t$ is defined as the argument of $\Psi$.

Let $f: [0,2\pi)\times [0,1]\to {\bf R}$ be defined by
$$f(t,s)=27 s^4-18 s^2+8s \cos t-1$$
and let $\omega(t)$ be the unique positive solution of the function $f(t,s)=0$ for each $t\in [0,2\pi)$. 

Hirakawa determined in \cite[2004]{Hir} the periodicity of $\phi:{\mathbb E}^2\to \tilde M^2(4\rho)$ as follows.

\begin{theorem} Let $\phi:{\mathbb E}^2\to CH^2(-4\rho)$ be an isometric immersion with nonzero parallel mean curvature vector determined by $b$ and $t$. Then

\vskip.04in
{\rm (1)} $\phi$ is doubly periodic if and only if $b$ and $t$ satisfy $\omega(t)<\sqrt{1+\rho/2b^2}$;

\vskip.04in
{\rm (2)} $\phi$ is singly periodic if and only if $b$ and $t$ satisfy $\sqrt{1+\rho/2b^2}\leq \omega(t)$ and $(t,\sqrt{1+\rho/2b^2})\ne (0,\omega(0))$;

\vskip.04in
{\rm (3)} $\phi$ is an imbedding if and only if $b$ and $t$ satisfy $(t,\sqrt{1+\rho/2b^2})=(0,\omega(0))$.

\end{theorem}

\section{Biharmonic surfaces  with parallel mean curvature vector.}

The study of biharmonic maps between two Riemannian manifolds, as a generalization
of harmonic maps, was suggested by  Eells and  Sampson in \cite[1965]{ES}.
Biharmonic maps $\psi:(M,g)\to (N,h)$ between Riemannian manifolds are critical
points of the bienergy functional:
$$E_2(\psi)=\frac{1}{2} \int_M |\tau(\psi)|^2 v_g,$$
where $\tau(\psi)={\rm trace}(\nabla d\psi)$ is the tension field of $\psi$. A biharmonic map is call proper if it is non-harmonic.

The  Euler-Lagrange equation for the biharmonic integral was computed by Jiang in  \cite[1986]{Jiang} which reads
$$\tau_2(\psi)=-\Delta \tau(\psi) -{\rm trace} R^N(d\psi,\tau(\psi))d\psi=0.$$

In another context,  in 1984 the author began the study of {\it biharmonic immersions} $\phi:M\to \mathbb E^m_s$ of a (pseudo) Riemannian manifold into (pseudo) Euclidean spaces defined by the property: $\Delta^2 \phi=0$. 

In Euclidean spaces this notion of biharmonic submanifolds introduced by the author agrees with the notion of biharmonic maps defined by Eells and Samson.

For biharmonic submanifolds in a pseudo-Euclidean space $\mathbb E^m_s$ with index $s$, the author derived a characterization of biharmonic submanifolds in terms of PDE system (see, for instance, \cite{c1984,ci1,ci2}). In particular, for biharmonic submanifolds in a Euclidean space, he proved the following.

\begin{proposition} \label{P:9.1}  An $n$-dimensional submanifold $M$ of $\mathbb E^m$ is biharmonic if and only if the following two
differential equations hold:
\begin{align}&\label{9.1} \Delta^{D} H + \sum_{i=1}^{n}h(A_{H}e_{i},e_{i})=0,
\\ & \label{9.2} n\nabla\<H,H\>+4\,{\rm trace}\,A_{DH}=0,\end{align}
where $\Delta^D$ is the Laplace operator associated with the normal connection $D$ and $\{e_1,\ldots,e_n\}$ is an orthonormal frame of $M$.
\end{proposition}

The following result  follows very easily from Proposition \ref{P:9.1}.

\begin{theorem} \label{P:9.1}  Let $M$ be $n$-dimensional submanifold of $\mathbb E^m$ with parallel mean curvature vector. Then $M$ is biharmonic if and only if $M$ is minimal.\end{theorem}

This follows from the fact that if $M$ has parallel mean curvature vector,  condition \eqref{9.2} holds  trivially. Moreover, condition \eqref{9.1} reduces to $ \sum_{i=1}^{n}h(A_{H}e_{i},e_{i})=0$, which implies that $H=0$, i.e., $M$ is minimal.

The author made the following conjecture in 1985  (see, for instance,  \cite[page 221]{report} and \cite[1992]{D}).
\vskip.06in

{\bf Chen's Conjecture.} {\it Every biharmonic submanifold  of a Euclidean space is a   minimal submanifold.}
\vskip.06in

During the last decade important progress has been made in the study of both the
geometry and the analytic properties of biharmonic maps and biharmonic submanifolds. In differential geometry,
a special attention has been payed to biharmonic submanifolds related with this conjecture. Very recently, Ou provided  in \cite[2010]{Ou} an equivalent formulation of this conjecture on biharmonic hypersurfaces by using biharmonic graph equation.

In the last decade, several differential geometers also investigated the following.
\vskip.06in

{\bf Generalized Chen's Conjecture.}  {\it Biharmonic submanifolds of a non positive sectional curvature manifold are minimal.}
\vskip.06in

This conjecture encouraged the study of proper biharmonic submanifolds in spheres and other curved spaces as well.

For biharmonic surfaces with parallel mean curvature vector in spheres,  Balmus, Montaldo and Oniciuc  proved in \cite[2008]{BMO} the following result.

\begin{theorem} Let $M$ be a proper biharmonic surface with parallel mean curvature
vector in $S^n(1)$. Then $M$  is minimal in a totally umbilical $S^{n-1}\big(\frac{1}{\sqrt{2}}\big)\subset S^n(1)$.\end{theorem}

\section{Submanifolds with parallel mean curvature vector satisfying additional conditions.}

It is a classical result of Liebmann that the only closed convex surfaces in Euclidean 3-space $\mathbb E^3$ having constant mean curvature are  round spheres. In \cite[1973]{Smyth}  Smyth  extended Liebmann's result to the following.

\begin{theorem} Let $M$ be a compact $n$-dimensional submanifold with nonnegative sectional curvature in Euclidean $m$-space. If  $M$ has parallel mean curvature vector, then $M$ is a product submanifold $M_1\times\cdots\times M_k$, where each $M_i$ is a minimal submanifold in a hypersphere of an affine subspace of $\mathbb E^m$.
\end{theorem}

Furthermore, Yano  and Ishihara \cite[1971]{YI} and Erbacher \cite[1972]{Er} extended Liebmann's result to

\begin{theorem} Let $M$ be an $n$-dimensional submanifold in Euclidean $m$-space with nonnegative sectional curvature. Suppose that the  mean curvature vector  is parallel in the normal bundle and the normal connection is flat. If $M$ is either compact or has constant scalar curvature, then   $M$ is the standard product immersion of the product  $S^{n_1}(r_1)\times\cdots\times S^{n_k}(r_k)$ of some spheres.
\end{theorem}

For complete submanifolds $M^n$ of dimension $\geq 3$ in Euclidean space,   Shen \cite[1985]{Sh} proved the following.

\begin{theorem} Let $M^n$ $(n\ge 3)$ be a  complete submanifold in the Euclidean space $\mathbb E^{m}$ with parallel mean curvature vector. If the squared mean curvature $H^2$ and the squared norm $S=||h||^2$ of the second fundamental form satisfies $$(n-1)S\leq n^2H^2,$$ then $M^n$ is an $n$-plane, an $n$-sphere $S^n$, or a circular cylinder $S^{n-1}\times \mathbb E^1$.
\end{theorem}

Shen's result extended a result of Chen-Okumura \cite[1973]{COk}.

G. Chen and X. Zou \cite[1995]{CZ} studied compact submanifolds of spheres with parallel mean curvature vector and proved the following.

\begin{theorem} Let $M$ be an $n$-dimensional compact submanifold with nonzero parallel mean curvature vector in the unit $(n+p)$-sphere. Then 
\vskip.04in

{\rm (1)} $M$ is totally geodesic, if one of the following two conditions hold:
\begin{equation}\begin{aligned} &S\leq \min\left\{\frac23 n,\frac{2n}{1+\sqrt{\frac n 2}}\right\},\quad
p\geq 2\quad \hbox{and}\quad n\ne 8;\\& \notag S\leq \min\left\{\frac n{2-\frac 1{p-1}},\frac{2n}{1+\sqrt{\frac n 2}}\right\},\quad
p\geq 1\quad\hbox{and}\quad (n,p)\ne (8,3);\end{aligned}\end{equation}
\vskip.04in

{\rm (2)} $M$ is totally umbilical, if $2\leq n\leq 7,\; p\geq 2,$ and $S\leq \frac23 n$.
\end{theorem}

In \cite[2001]{CN}, Cheng and Nonaka  proved the following.

\begin{theorem} \label{T:11.5} Let $M$ be an $n$-dimensional complete and connected submanifolds with parallel mean curvature vector in the $(n+p)$-dimensional Euclidean space $\mathbb E^{n+p}$.  If the squared norm of the second fundamental form $h$ satisfies $$S \leq \frac{n^2 H^2}{n-1},$$ then it is the totally geodesic Euclidean space $\mathbb E^n$, the totally umbilical sphere $S^n(c)$ or the generalized cylinder $S^{n-1}(c)\times \mathbb E^1$. 
\end{theorem}

Theorem \ref{T:11.5} result is a generalization of a classical result of Klotz and Osserman \cite[1966/7]{KO}.

\section{Lagrangian submanifolds with parallel mean curvature vector} 

A totally real submanifold $M$ in a Kaehler manifold $\tilde M$ is called {\it Lagrangian} if $\dim M=\dim_{\bf C}\tilde M$. 

The following result was obtained by Chen, Houh and Lue in \cite[1977]{CHL}.

\begin{theorem}\label{T:8.1}  Let $M$ be a compact Lagrangian submanifold  of ${\bf C}^n$. If $M$ has nonnegative sectional curvature and parallel mean curvature vector, then $M$ is a product submanifold $M_1\times\cdots\times M_r$, where  $M_j$ is a compact Lagrangian submanifold embedded in some  complex linear subspace {\bf C}$^{n_j}$ and is immersed as a minimal submanifold in a hypersphere of {\bf C}$^{n_j}$.
\end{theorem}

Theorem \ref{T:8.1} was extended to complete Lagrangian minimal submanifolds of $\hbox{\bf C}^n$ by  Urbano in  \cite[1995]{U1} and by Ki and Kim in \cite[1996]{KK}

 Ohnita \cite[1986]{On} and  Urbano \cite[1986]{U2} investigated Lagrangian submanifolds  with $K\geq 0$ in complex projective space and obtained the following.

\begin{theorem} \label{T:8.2} Let $M$ be a compact Lagrangian  submanifold of  $CP^n(4)$ with  parallel mean
curvature vector. If $M$ has nonnegative sectional curvature, then the second fundamental form of $M$ is parallel. \end{theorem} 

By applying Theorem \ref{T:8.2} and the classification of parallel submanifolds of $CP^n$ obtained by Naitoh in \cite[1983]{Nai}, it follows that
compact Lagrangian submanifolds of $CP^n(4)$ with parallel mean curvature vector and $K\geq 0$ must be the product submanifolds: $$T\times M_1\times\cdots\times M_k,$$ where $T$ is a flat torus with
$\dim M\geq k-1$ and each $M_i$ is one
of the following:
\begin{equation}\begin{aligned} RP^r(1)\to CP^r(4)\;\;(r\geq 2)\;\;(\hbox{totally geodesic}),\\
SU(r)/SO(r)\to CP^{(r-1)(r+2)/2}(4)\;\;(r\geq 3)\;\;(\hbox{minimal}),\\ SU(2r)/Sp(r)\to CP^{(r-1)(2r+1)}(4)\;\;(r\geq 3)\;\;(\hbox{minimal}),\\ SU(r)\to CP^{r^2-1}(4)\;\;(r\geq 3)\;\;(\hbox{minimal}),\\ E_6/F_4\to CP^{26}(4)\;\;(\hbox{minimal}).\\ \end{aligned}\end{equation}

For $3$-dimensional compact Lagrangian   submanifolds of complex space forms,
 Urbano proved  in \cite[1996]{U2} the following.

\begin{theorem}  If $M^3$ is a $3$-dimensional compact Lagrangian  submanifold of a complex space form $\tilde M^3(4c)$ with nonzero parallel mean curvature vector, then $M^3$ is flat and has parallel second fundamental form. \end{theorem}

\section{Black holes and  trapped surfaces.}

The idea of an object with gravity strong enough to prevent light from escaping was proposed in 1783 by  J. Michell (1724 - 1793), an amateur British astronomer. P.-S. Laplace (1749 - 1823), a French physicist  came to the same conclusion in 1795 independently. Black holes, as currently understood, are described by Einstein's general theory of relativity developed in  \cite[1916]{E}. 

Einstein's  theory of general relativity predicts that when a large enough amount of mass is present in a sufficiently small region of space, all paths through space are warped inwards towards the center of the volume, preventing all matter and radiation within it from escaping. Einstein's theory has important astrophysical applications. It points towards the existence of black holes. In addition, general relativity is the basis of current cosmological models of an expanding universe.

According to the American Astronomical Society,  every large galaxy has a super massive black hole  ($\sim 10^5$ - $10^9\, M_{sun}$) at its center.
The black hole's mass is proportional to the mass of the host galaxy,  suggesting that the two are linked very closely. 

Black holes can't be seen, because everything that falls into them, including light, is trapped. 
 But the swift motions of gas and stars near an otherwise invisible object allows astronomers to calculate that it's a black hole and even to estimate its mass. 

The concept of {\it trapped surfaces}, introduced by Penrose in \cite[1965]{Pe} plays an  important role in general relativity. In the theory of  cosmic black holes, if there is a massive source inside the surface, then close enough to a massive enough source, the outgoing light rays may also be converging; a {\it trapped surface}. Everything inside is trapped. Nothing can escape, not even light. It is believed that there will be a marginally trapped surface, separating the trapped surfaces from the untrapped ones, where  the outgoing light rays are instantaneously parallel. The surface of a {\it black hole} is located by the marginally trapped surface. 
In terms of the mean curvature vector, a  surface  is {\it marginally trapped} if and only if its mean curvature vector is light-like at each point on the surface. 

In the last few years, marginally trapped surfaces have been studied extensively from differential geometric view point.
For recent  survey on results of marginally trapped surfaces from this view point, see \cite[2009]{IEJG} and \cite[2009]{Tamkang}.

\section{Marginally trapped surfaces with parallel mean curvature 4D Lorentzian space forms.}

Marginally trapped surfaces with parallel mean curvature vector in 4-dimensional Lorentzian space form were  classified by the author and Van der Veken in \cite[2010]{cv3}. More precisely, they obtained the following three results.
 
\begin{theorem}\label{T:12.1} Let  $M$ be a marginally trapped surface  with parallel mean curvature vector in the Minkowski space-time $\mathbb E^4_1$. Then, with respect to  suitable Minkowskian coordinates $(t,x_2,x_3,x_4)$ on $\mathbb E^4_1$, $M$ is  an open part of one of the following six types of surfaces:
\vskip.05in

 {\rm (1)}  a flat parallel biharmonic surface given by
$$\text{\small$\frac{1}{2}$}\Big((1-b)u^2+(1+b)v^2,(1-b)u^2+(1+b)v^2,2u,2v\Big),\;
b\in \bf R.$$ \vskip.05in

 {\rm (2)}  a flat parallel surface given by $$a\big(\cosh
u, \sinh u,\cos v,\sin v\big),\; a>0.$$ \vskip.05in

 {\rm (3)}  a  flat surface given by $(f(u,v),u,v,f(u,v)),$ where $f$ is a function on $M$ such that $\Delta f=c$ for some nonzero real number $c$, where $\Delta$ is the Laplacian of $M$.
 \vskip.05in

 {\rm (4)}  a   flat  surface lying in the light cone $\mathcal{LC}$.
\vskip.05in

 {\rm (5)}  a non-parallel surface lying in the de Sitter space-time $S^3_1(r^2)$ for some $r>0$ such that the  mean curvature
vector $H'$ of $M$ in $S^3_1(r^2)$ satisfies $\<H',H'\>=-r^2$.
\vskip.05in

 {\rm (6)}  a non-parallel surface lying in the hyperbolic space $H^3(-r^2)$ for some $ r>0$ such that the  mean curvature
vector $H'$ of $M$ in $H^3(-r^2)$ satisfies $\<H',H'\>=r^2$.
\vskip.05in

Conversely, all surfaces of types $(1)$--$(6)$  above give rise to marginally trapped surfaces with parallel mean curvature vector in
$\mathbb E^4_1$.
\end{theorem}

For marginally trapped surfaces \with in the de Sitter 4-space, we have the following.

\begin{theorem}\label{T:12.2} Let  $M$ be a marginally trapped  surface  with parallel mean curvature vector in the de Sitter space-time $S^4_1(1)\subset \mathbb E^5_1$.  Then $M$ is congruent to an open part of one of the following eight types of surfaces:
\vskip.05in

 {\rm (i)} a  parallel surface of curvature one given by
$$\big(1,\sin u,\cos u \cos v,\cos u\sin v,1\big).$$

\vskip.05in  {\rm (ii)}   a flat parallel surface defined
by $$\frac{1}{2}\big(2u^2-1,2u^2-2,2u,\sin 2v,\cos 2v\big).$$

 {\rm (iii)}   a flat parallel surface defined by
\begin{equation}\begin{aligned}\notag
&\text{\small$\Bigg(\frac{b}{\sqrt{4-b^2}},\frac{\cos (\sqrt{2-b}\,u)}{\sqrt{2-b}},\frac{\sin (\sqrt{2-b}\,u)}{\sqrt{2-b}},\frac{\cos
(\sqrt{2+b}\,v)}{\sqrt{2+b}},\frac{\sin (\sqrt{2+b}\,v)}{\sqrt{2+b}} \Bigg)$},\,  \end{aligned}\end{equation}
with $|b|<2.$

 {\rm (iv)}    a flat parallel surface defined by
\begin{equation}\begin{aligned}\notag
&\text{\small$\Bigg(\frac{\cosh (\sqrt{b-2}\,u)}{\sqrt{b-2}},\frac{\sinh (\sqrt{b-2}\,u)}{\sqrt{b-2}},\frac{\cos
(\sqrt{2+b}\,v)}{\sqrt{2+b}},\frac{\sin
(\sqrt{2+b}\,v)}{\sqrt{2+b}}, \frac{b}{\sqrt{b^2-4}} \Bigg)$},\,
\end{aligned}\end{equation} with $ b>2.$

\vskip.05in {\rm (v)}  a surface of constant curvature one immersed in $S^4_1(1)\subset \mathbb E^5_1$ given by
$$(f,\cos u,\sin u\cos v,\sin u \sin v,f),$$
where $f$ is a function satisfies $\Delta f=k$ for some nonzero real number $k$.
 \vskip.05in

 {\rm (vi)}   a non-parallel surface of
curvature one in $S^4_1(1)$ which lies in  $$\mathcal {LC}_1:=\{({\bf
y},1)\in \mathbb E^5_1: \<{\bf y},{\bf y}\>=0,\ {\bf y}\in \mathbb
E^4_1\} \subset S^4_1(1);$$

\vskip.05in
 {\rm (vii)}  a non-parallel surface  of
$S^4_1(1)$ which lies  in $S^4_1(1)\cap S^4_1(c_0,r^2)$ with $c_0\ne
0$ and $r>0$ such that the mean curvature vector $H'$ of $M$ in
$S^4_1(1)\cap S^4_1(c_0,r^2)$ satisfies $\<H',H'\>=1-r^2$.
\vskip.05in

 {\rm (viii)}   a non-parallel surface of $S^4_1(1)$ which
lies in $S^4_1(1)\cap H^4(c_0,-r^2)$ with $c_0\ne 0$ and $r>0$ such
that the mean curvature vector $H'$ of $M$ in $S^4_1(1)\cap
H^4(c_0,-r^2)$ satisfies $\<H',H'\>=1+r^2$.

Conversely, all surfaces of types ${\rm (i)}$--${\rm (viii)}$  above
give rise to marginally trapped surfaces with parallel mean
curvature vector in $S^4_1(1)$.

\end{theorem}

For marginally trapped surfaces \with in the anti-de Sitter 4-space, we have the following.

\begin{theorem}\label{T:12.3}Let  $M$ be a marginally trapped  surface  with parallel mean curvature vector in the anti de Sitter space-time $H^4_1(-1)\subset \mathbb E^5_2$.  Then, $M$ is congruent to one of the following eight types of surfaces:

 {\rm (a)}  a curvature $-1$ parallel surface  given by
$$\big(1,\cosh u \cosh v,\sinh u,\cosh u\sinh v,  1\big).$$

\vskip.05in
 {\rm (b)}   a flat parallel surface defined by $$\frac{1}{2}\big(2u^2+2,\cosh 2v, 2u,\sinh
2v,2u^2+1\big).$$

 {\rm (c)}   a flat parallel surface defined by
\begin{equation}\begin{aligned}\notag &  \text{\small$\Bigg(\frac{\cosh (\sqrt{2\hskip-.02in -\hskip-.02in b}\,u)}{\sqrt{2-b}},\frac{\cosh (\sqrt{2\hskip-.02in +\hskip-.02in b}\,v)}{\sqrt{2+b}}, \frac{\sinh (\sqrt{2\hskip-.02in -\hskip-.02in b}\,u)}{\sqrt{2-b}},\frac{\sinh (\sqrt{2\hskip-.02in +\hskip-.02in b}\,v)}{\sqrt{2+b}},\frac{b}{\sqrt{4\hskip-.02in -\hskip-.02in
b^2}}, \Bigg)$},\,  \end{aligned}\end{equation} with $|b|<2.$

 {\rm (d)}    a flat parallel surface defined by \begin{equation}\begin{aligned}\notag
&\text{\small$\Bigg(\frac{b}{\sqrt{b^2-4}},\frac{\cosh (\sqrt{b+2}\,v)}{\sqrt{b+2}},\frac{\sinh (\sqrt{b+2}\,v)}{\sqrt{b+2}}, \frac{\cos
(\sqrt{b-2}\,u)}{\sqrt{b-2}},\frac{\sin (\sqrt{b-2}\,u)}{\sqrt{b-2}} \Bigg)$},\, \end{aligned}\end{equation} with $ b>2.$

\vskip.05in
 $(e)$  a surface of constant curvature $-1$ immersed in $H^4_1(-1)\subset \mathbb E^5_2$ given by
$$(f,\cosh u,\sinh u\cos v,\sinh u \sin v,f),$$
where $f$ is a function satisfies $\Delta f=k$ for some nonzero real number $k$.
\vskip.05in

 $(f)$  a non-parallel surface of $H^4_1(-1)$ with  curvature $-1$, lying in $$\mathcal
{LC}_2:=\{(1,{\bf y})\in \mathbb E^5_2: \<{\bf y},{\bf y}\>=0,\ {\bf y}\in \mathbb E^4_1\} \subset H^4_1(-1).$$
\vskip.05in

 $(g)$ a non-parallel surface lying in $H^4_1(-1)\cap S^4_2(c_0,r^2)$ with $c_0\ne 0$ and $r>0$ such that the mean
curvature vector $H'$  in $H^4_1(-1)\cap S^4_2(c_0,r^2)$ satisfies $\<H',H'\>=-r^2-1$;.\vskip.05in

 $(h)$  a non-parallel surface lying  in $H^4_1(-1)\cap H_1^4(c_0,-r^2)$ with $c_0\ne 0$ and $r>0$ such that mean curvature
vector $H'$   in $H^4_1(-1)\cap H_1^4(c_0,-r^2)$ satisfies $\<H',H'\>=r^2-1$.

Conversely, all surfaces of types $(a)$--$(h)$  above give rise to marginally trapped surfaces with parallel mean curvature vector in $H^4_1(-1)$.
\end{theorem}

\section{Classification of space-like surfaces with parallel mean curvature vector in $\mathbb E^m_s$.}

Space-like surfaces with parallel mean curvature vector in Lorentzian space forms $L^4(c)$ were classified by the author and Van der Veken in \cite{cv}.
Space-like surfaces with parallel mean curvature vector in a pseudo-Euclidean space with arbitrary codimension and arbitrary index were completely classified by the author in \cite[2009]{c6}. More precisely, he proved the following.

\begin{theorem}\label{T:15.1} Let $\phi: M\to \mathbb E^m_s$ be an isometric immersion of a space-like surface into a pseudo-Euclidean $m$-space $\mathbb E^m_s$ with index $ s\geq 1$. Then $M$ has parallel mean curvature vector if and only if it is congruent to one of the following sixteen types of surfaces:

\begin{enumerate}

\vskip.04in
\item[(1)] a minimal surface of $\, \mathbb E^m_s$;

\vskip.04in
\item[(2)]  a minimal surface of a pseudo-Riemannian $(m-1)$-sphere $S^{m-1}_s(r^2)\subset \mathbb E^m_s$ with $ r>0$;

\vskip.04in
\item[(3)]  a minimal surface of a pseudo-hyperbolic $(m-1)$-space $H^{m-1}_{s-1}(-r^2)\subset \mathbb E^m_s$ with $ r>0$;

\vskip.04in
\item[(4)]   a CMC surface of a  Euclidean $3$-space $\mathbb E^3\subset \mathbb E^m_s$;

\item[(5)] a CMC surface of a  Minkowski 3-space $\mathbb E^3_1\subset \mathbb E^m_s$;

\vskip.04in
\item[(6)]   a CMC surface of a $3$-sphere  $S^3(r^2),\, r>0$, where $S^3(r^2)$ is an ordinary hypersphere of a  Euclidean 4-space $ \mathbb E^4\subset \mathbb E^m_s\, (s\leq m-4)$;.

\vskip.04in
\item[(7)] a CMC surface  of a de Sitter 3-space $S^3_1(r^2)\subset \mathbb E^4_1\subset \mathbb E^m_s$;

\vskip.04in
\item[(8)]  a CMC surface  of a  hyperbolic 3-space $H^3(-r^2)\subset \mathbb E^4_1\subset \mathbb E^m_s$;

\vskip.04in
\item[(9)]  a CMC surface of an anti-de Sitter 3-space $H^3_1(-r^2) \subset \mathbb E^4_2\subset \mathbb E^m_s$, $r>0$; 

\vskip.04in
\item[(10)]  a  surface with constant Gauss curvature in the light cone $\mathcal{LC}\subset \mathbb E^4_1\subset \mathbb E^m_s$;

\vskip.04in
\item[(11)] a flat  surface given by $(a\cosh u, a\sinh u,a\cos v,a\sin v)\in \mathbb E^4_1\subset \mathbb E^m_s$,  $a>0;$

\vskip.04in
\item[(12)] a marginally trapped surface given by $\phi=(f,\psi,f)$, where $f$ is a function on $M$ satisfying $\Delta f=b$ for some real number $b\ne 0$ and $\psi:M\to \mathbb E^{m-2}_{s-1}$ is a minimal immersion;

\vskip.04in
\item[(13)]   a surface given by  $\phi=(f,\psi,f),$  where  $f$  is a function satisfying $\Delta f=k$ for some  $k\in {\bf R}$ and $\psi:M \to \mathbb E^3\subset \mathbb E^{m-2}_{s-1}$ is a CMC surface.

\vskip.04in
\item[(14)]  
 a surface given by  $\phi=(f,\psi,f),$  where  $f$  is a function satisfying $\Delta f=k$ for some  $k\in {\bf R}$ and $\psi:M \to \mathbb E^3_1\subset \mathbb E^{m-2}_{s-1}$ is a CMC surface.

\vskip.04in
\item[(15)]   a non-parallel surface lying in the de Sitter
space-time $S^3_1(r^2)\subset \mathbb E^4_1\subset \mathbb E^m_s$ for some $r>0$ such that the  mean curvature
vector $H'$ of $M$ in $S^3_1(r^2)$ satisfies $\<H',H'\>=-r^2$;

\vskip.04in
\item[(16)] a non-parallel surface lying in the hyperbolic
space $H^3(-r^2)\subset \mathbb E^4_1\subset \mathbb E^m_s$ for some $ r>0$ such that the  mean curvature vector
$H'$ of $M$ in $H^3(-r^2)$ satisfies $\<H',H'\>=r^2$.
\end{enumerate}
In the statements above, $\mathbb E^3_i,\, \mathbb E^4_i$ and $\mathbb E^5_i$ $(i=0,1,2)$ are imbedded in $\mathbb E^m_s$ as totally geodesic submanifolds in standard way.
\end{theorem} 

    As an immediate by-product of Theorem \ref{T:15.1} we get the following complete classification of space-like surfaces with parallel mean curvature vector in Minkowski spaces of arbitrary dimension.
    
\begin{theorem}\label{T:15.2} Let $\phi: M\to \mathbb E^m_1$ be an isometric immersion of a space-like surface into a Minkowski $m$-space $\mathbb E^m_1$. Then $M$ has parallel mean curvature vector if and only if it is congruent to one of the following fourteen types of surfaces:

\begin{enumerate}

\vskip.04in
\item[(1)] a minimal surface of $\, \mathbb E^m_1$;

\vskip.04in
\item[(2)]  a minimal surface of a de Sitter $(m-1)$-space $S^{m-1}_1(r^2)$ with $ r>0$;

\vskip.04in
\item[(3)]  a minimal surface of a hyperbolic $(m-1)$-space $H^{m-1}(-r^2)$ with $ r>0$;

\vskip.04in
\item[(4)]   a CMC surface of a  Euclidean $3$-space $\mathbb E^3\subset \mathbb E^m_1$;

\item[(5)] a CMC surface of a  Minkowski 3-space $\mathbb E^3_1\subset \mathbb E^m_1$;

\vskip.04in
\item[(6)]   a CMC surface of a $3$-sphere  $S^3(r^2),\, r>0$, where $S^3(r^2)$ is an ordinary hypersphere of a  Euclidean 4-space $ \mathbb E^4\subset \mathbb E^m_1$;.

\vskip.04in
\item[(7)] a CMC surface  of a de Sitter 3-space $S^3_1(r^2)\subset \mathbb E^4_1\subset \mathbb E^m_1$;

\vskip.04in
\item[(8)]  a CMC surface  of a  hyperbolic 3-space $H^3(-r^2)\subset \mathbb E^4_1\subset \mathbb E^m_1$;

\vskip.04in
\item[(9)] a flat  surface given by $(a\cosh u, a\sinh u,a\cos v,a\sin v)\in \mathbb E^4_1\subset \mathbb E^m_1$,  $a>0;$

\vskip.04in
\item[(10)]  a  surface with constant Gauss curvature in the light cone $\mathcal{LC}\subset \mathbb E^4_1\subset \mathbb E^m_1$;

\vskip.04in
\item[(11)] a marginally trapped surface given by $\phi=(f,\psi,f)$, where $f$ is a function on $M$ satisfying $\Delta f=b$ for some real number $b\ne 0$ and $\psi:M\to \mathbb E^{m-2}$ is a minimal immersion;

\vskip.04in
\item[(12)]  
  a surface given by  $\phi=(f,0,\ldots,0,\psi,f),$  where  $f$  is a function satisfying $\Delta f=k \, (k\in {\bf R})$ and $\psi:M \to \mathbb E^3$ is a CMC surface.
  
\vskip.04in
\item[(13)]   a non-parallel surface lying in the de Sitter
space-time $S^3_1(r^2)\subset \mathbb E^4_1\subset \mathbb E^m_1$ for some $r>0$ such that the  mean curvature
vector $H'$ of $M$ in $S^3_1(r^2)$ satisfies $\<H',H'\>=-r^2$;

\vskip.04in
\item[(14)] a non-parallel surface lying in the hyperbolic
space $H^3(-r^2)\subset \mathbb E^4_1\subset \mathbb E^m_1$ for some $ r>0$ such that the  mean curvature vector
$H'$ of $M$ in $H^3(-r^2)$ satisfies $\<H',H'\>=r^2$.

\end{enumerate}
In the statements above, $\mathbb E^3,\, \mathbb E^3_1,\, \mathbb E^4,\, \mathbb E^4_1$ and $\mathbb E^5_1$ are imbedded in $\mathbb E^m_1s$ as totally geodesic submanifolds in standard way.
\end{theorem} 

\section{Space-like surfaces with parallel mean curvature vector in $S^m_s(1)$.}

Space-like surfaces \with in pseudo-Riemannian spheres with arbitrary codimension and  index were completely classified by the author in \cite[2009]{cent}.

\begin{theorem}\label{T:16.1} Let $\psi:M\to S^m_s(1)$ be an isometric immersion of a space-like surface into a pseudo-Riemannian $m$-sphere $S^m_s(1)$ with index $s\geq 1$. If $M$ has parallel mean curvature vector, then $M$ is locally congruent to an open portion of one of the following fifteen types of surfaces:

\begin{enumerate}

\vskip.03in
\item[(1)] a minimal surface of $S^m_s(1)$;

\vskip.03in
\item[(2)] a surface given by $(f,\phi,f)$, where $f$ is function satisfying $\Delta f=2f-r$ with $0\ne  r\in {\bf R}$ and $\phi:M\to S^{m-2}_{s-1}(1)$  is a minimal surface;

\vskip.03in
\item[(3)]  a minimal surface of  a pseudo-Riemannian sphere $S^{m-1}_s$  imbedded  as \begin{align}\notag &\hskip.4in \left\{\Big(y,\text{\small$\frac{\alpha}{\sqrt{1\!+\!\alpha^{2}}}$}\Big)\in \mathbb E^{m+1}_s:y\in \mathbb E^m_s \; and \; \<y,y\>=\text{\small$ \frac{1}{1+\alpha^2}$}\right\},\;  \alpha>0;\end{align}

\vskip.03in
\item[(4)]  a minimal surface of a pseudo-Riemannian sphere $S^{m-1}_{s-1}$ imbedded  as \begin{align}\notag &\hskip.4in \left\{\Big(\text{\small$\frac{\alpha}{\sqrt{1-\alpha^2}}$},y\Big)\in \mathbb E^{m+1}_s: y\in \mathbb E^m_{s-1}\; and \;  \<y,y\>=\text{\small$ \frac{1}{1-\alpha^2}$}\right\}\end{align} with $ \alpha\in (0,1)$;

\vskip.03in
\item[(5)]  a minimal surface of a pseudo-hyperbolic space $H^{m-1}_{s-1}$ imbedded as \begin{align}\notag &\hskip.4in \left\{\Big(y,\text{\small$\frac{\alpha}{\sqrt{1\!+\!\alpha^{2}}}$}\Big)\in \mathbb E^{m+1}_s:y\in \mathbb E^m_{s}\; and \; \<y,y\>=\text{\small$ \frac{1}{1-\alpha^2}$}\right\},\;\; \alpha >1;\end{align}

\vskip.03in
\item[(6)] A surface with parallel mean curvature vector  in a 4-sphere $S^4$ imbedded in $S^m_s(1)$ as $$\{(b,y)\in \mathbb E^{m+1}_s:y\in \mathbb E^5, \<b,b\>=1-r^2, \<y,y\>=r^2>0\}$$ with $b\in \mathbb E^{m-4}_s$; 

\item[(7)] a surface with parallel mean curvature vector of a de Sitter 4-space $S^4_1$, which is imbedded in  $S^m_s(1)$ as 
$$ \left\{(b,y)\in \mathbb E^{m+1}_s:y\in \mathbb E^5_{1}, \<y,y\>=1-\<b,b\> \right\}$$  for some vector  $b\in \mathbb E^{m-4}_{s-1}$ with $\<b,b\><1$;

\vskip.03in
\item[(8)] a surface with parallel mean curvature vector lying in a pseudo-Riemannian 4-sphere $S^4_2$ imbedded as 
$$ \{(b,y)\!\in\! \mathbb E^{m+1}_s:y\in \mathbb E^5_{2}, \<y,y\> =1-\<b,b\>\}$$ for some 
 $b\in \mathbb E^{m-4}_{s-2}$ with $\<b,b\><1$;

\vskip.03in
\item[(9)] a surface with parallel mean curvature vector in a hyperbolic 4-space imbedded in  $S^m_s(1)\!$ as 
$$\left\{(b,y)\in  \mathbb E^{m+1}_s: y\in \mathbb E^5_{1}, \<y,y\>=1-\<b,b\>\right\}$$ for some $b\in \mathbb E^{m-4}_{s-1}$ with $\<b,b\>>1$;

\vskip.03in
\item[(10)] a surface with parallel mean curvature vector of an anti-de Sitter  4-space $H^4_1$  imbedded  in  $S^m_s(1)$ as 
$$ \left\{(b,y)\in   \mathbb E^{m+1}_s:y\in \mathbb E^5_{2}, \<y,y\>=1-\<b,b\>\right\}$$ for  $b\in \mathbb E^{m-4}_{s-2}$ with $\<b,b\>>1$;

\vskip.03in
\item[(11)] a surface of constant Gauss curvature lying in the light cone ${\mathcal LC}^3_1\subset {\bf R}^4_{1,0}$ and ${\mathcal LC}^3_1$ is  imbedded in  $S^m_s(1)$ as 
$$\left\{(b,y)\in \mathbb E^{m+1}_s:  y\in \mathbb E^4_1, \<y,y\>=0\right\},$$
where $b$ is a unit space-like vector in $\mathbb E^{m-3}_{s-1}$;

\vskip.03in
\item[(12)] a surface  defined by $(f,b,\psi,f): M\to \mathbb E^{m+1}_s \, (s\leq m-4)$, where  $b\in \mathbb E^{m-5}_{s-1}$ satisfies $\<b,b\><1$,  $\psi:M\to S^3(1/(1-\<b,b\>))\subset \mathbb E^4$ has constant mean curvature and $f$ is a function satisfying $\Delta f=2f+2k$ for some $k\in {\bf R}$;

\vskip.03in
\item[(13)] a surface  defined by $(f,b,\psi,f): M\to \mathbb E^{m+1}_s \, (s\in [2,m-3])$, where  $b\in \mathbb E^{m-5}_{s-2}$ satisfies $\<b,b\><1$,  $\psi:M\to S^3_1(1/(1-\<b,b\>))\subset \mathbb E^4_1$ has constant mean curvature and $f$  satisfies $\Delta f=2f+2k$ for some $k\in {\bf R}$;

\vskip.03in
\item[(14)] a surface  defined by $(f,b,\psi,f): M\to \mathbb E^{m+1}_s$, where  $b\in \mathbb E^{m-5}_{s-2}$ satisfies $\<b,b\>>1$,  $\psi:M\to H^3(1/(1-\<b,b\>))\subset \mathbb E^4_1$ has constant mean curvature and $f$  satisfies $\Delta f=2f+2k$ for some $k\in {\bf R}$;

\vskip.03in
\item[(15)] a surface  defined by $(f,b,\psi,f): M\to \mathbb E^{m+1}_s \, (2\leq s\leq m-3)$, where  $b$ is a unit space-like vector in $ \mathbb E^{m-5}_{s-2}$,  $\psi:M\to {\mathcal LC}^3_1\subset \mathbb E^4_1$ has constant mean curvature and $f$  satisfies $\Delta f=2f+2k$ for some $k\in {\bf R}$.

\end{enumerate}

\end{theorem} 

\section{Space-like surfaces with parallel mean curvature vector in $H^m_s(-1)$.}

 In \cite[2009]{cent}, the author also classified space-like surfaces \with in pseudo-hyperbolic spaces with arbitrary codimension and  index.

\begin{theorem}\label{T:17.1} Let $\psi:M\to H^m_s(-1)$ be an isometric immersion of a spatial surface $M$ into the pseudo-hyperbolic $m$-space $H^m_s(-1)$. If $M$ has parallel mean curvature vector, then $M$ is  congruent to an open portion of one of the following sixteen types of surfaces:
\begin{enumerate}

\vskip.03in
\item[(1)] a minimal surface of $H^m_s(-1)$;

\vskip.03in
\item[(2)] a surface given by $(f,\phi,f)$, where $f$ is a function satisfying $\Delta f=r-2f$ with $0\ne r\in {\bf R}$;
 and $\phi:M\to H^{m-2}_{s-1}(-1)$  is a minimal surface;

\vskip.03in
\item[(3)]  a minimal surface of  a pseudo-Riemannian sphere $S^{m-1}_s$  imbedded  as \begin{align}\notag &\hskip.4in \left\{\Big(\text{\small$\frac{\alpha}{\sqrt{\alpha^{2}-1}}$},y\Big)\in \mathbb E^{m+1}_s:y\in \mathbb E^m_{s+1} \; and \; \<y,y\>=\text{\small$ \frac{1}{\alpha^2-1}$}\right\},\;  \alpha>1;\end{align}

\vskip.03in
\item[(4)]  a minimal surface of a pseudo-hyperbolic space $H^{m-1}_{s}$ imbedded  as
\begin{align}\notag &\hskip.4in\left\{\(y,\text{\small$\frac{\alpha}{\sqrt{1-\alpha^{2}}}$}\) \in \mathbb E^{m+1}_{s+1}:y\in \mathbb E^m_{s+1},\<y,y\>=\text{\small$\frac{-1}{1-\alpha^2}$}\right\},\, 0<\alpha<1; \end{align}
\vskip.03in

\item[(5)]  a minimal surface of a pseudo-hyperbolic space $H^{m-1}_{s-1}$ imbedded as 
\begin{align}\notag&\hskip.4in  \left\{\(\text{\small$\frac{\alpha}{\sqrt{1+\alpha^{2}}}$},y\)\in \mathbb E^{m+1}_{s+1}:y\in \mathbb E^m_{s}, \<y,y\>=\text{\small$\frac{-1}{1+\alpha^2}$}\right\},\;\; \alpha\ne 0;  \end{align}

\vskip.03in
\item[(6)] A surface with parallel mean curvature vector in a 4-sphere $S^4$ imbedded in $H^m_s(-1)$ as $$\{(c,y)\in \mathbb E^{m+1}_{s+1}: y\in \mathbb E^5,\<c,c\>=-1-r^2, \<y,y\>=r^2>0\},$$ with  $c\in \mathbb E^{m-4}_{s+1}$;

\item[(7)] a surface with parallel mean curvature vector of a de Sitter 4-space $S^4_1$ imbedded in  $H^m_s(-1)$ as 
$$ \left\{(b,y)\in \mathbb E^{m+1}_{s+1}:y\in \mathbb E^5_{1}, \<y,y\>=-1-\<b,b\> \right\}$$  for some vector  $b\in \mathbb E^{m-4}_{s}$ with $\<b,b\><-1$;

\vskip.03in
\item[(8)] a surface with parallel mean curvature vector lying in a pseudo-Riemannian 4-sphere $S^4_2$ imbedded as 
$$\{(b,y)\!\in\! \mathbb E^{m+1}_{s+1}: y\in \mathbb E^5_{2}, \<y,y\> =-1-\<b,b\>\}$$
 for some $b\in \mathbb E^{m-4}_{s-1}$ with $\<b,b\><-1$;

\vskip.03in
\item[(9)] a surface with parallel mean curvature vector of a hyperbolic  4-space $H^4$  imbedded in  $H^m_s(-1)$ as 
$$\left\{(b,y)\in \mathbb E^{m+1}_{s+1}: y\in \mathbb E^5_{1}, \<y,y\>=-1-\<b,b\>\right\}$$ for some $b\in \mathbb E^{m-4}_{s}$ with $\<b,b\>>-1$;

\vskip.03in
\item[(10)] a surface with parallel mean curvature vector of an anti-de Sitter  4-space $H^4_1$  imbedded  in  $H^m_s(-1)$ as 
$$ \left\{(b,y)\in   \mathbb E^{m+1}_{s+1}:y\in \mathbb E^5_{2}, \<y,y\>=-1-\<b,b\>\right\}$$ for some $b\in \mathbb E^{m-4}_{s-1}$ with $\<b,b\>>-1$;

\vskip.03in
\item[(11)] a surface with parallel mean curvature vector of a pseudo-hyperbolic  4-space $H^4_2$ imbedded  in  $H^m_s(-1)$ as 
$$ \left\{(b,y)\in   \mathbb E^{m+1}_{s+1}:y\in \mathbb E^5_{3},  \<y,y\>=-1-\<b,b\>\right\}$$ for some $b\in \mathbb E^{m-4}_{s-2}$ with $\<b,b\>>-1$;

\vskip.03in
\item[(12)] a surface of constant Gauss curvature lying in the light cone ${\mathcal LC}^3_1\subset {\bf R}^4_{1,0}$ and ${\mathcal LC}^3_1$ is  imbedded in  $H^m_s(-1)$ as  $$\left\{(b,y)\in \mathbb E^{m+1}_{s+1}: y\in \mathbb E^4_1,  \<y,y\>=0\right\}$$ with $b$ being a unit time-like vector in $\mathbb E^{m-3}_{s}$;

\vskip.03in
\item[(13)] a surface  defined by $(f,b,\psi,f): M\to \mathbb E^{m+1}_{s+1} \, (s\geq 2)$, where  $b\in \mathbb E^{m-5}_{s-1}$ satisfies $\<b,b\><-1$,  $\psi:M\to S^3_1(-1/(1+\<b,b\>))\subset \mathbb E^4_1$ has constant mean curvature and $f$  satisfies $\Delta f=2k-2f$ for some $k\in {\bf R}$;

\vskip.03in
\item[(14)] a surface  defined by $(f,b,\psi,f): M\to \mathbb E^{m+1}_{s+1}$, where  $b\in \mathbb E^{m-5}_{s-1}$ satisfies $\<b,b\>>-1$,  $\psi:M\to H^3(-1/(1+\<b,b\>))\subset \mathbb E^4_1$ has constant mean curvature and $f$  satisfies $\Delta f=2k-2f$ for some $k\in {\bf R}$;

\vskip.03in
\item[(15)] a surface  defined by $(f,b,\psi,f): M\to \mathbb E^{m+1}_{s+1}$, where  $b\in \mathbb E^{m-5}_{s-2}$ satisfies $\<b,b\>>-1$,  $\psi:M\to H^3_1(-1/(1+\<b,b\>))\subset \mathbb E^4_2$ has constant mean curvature and $f$  satisfies $\Delta f=2k-2f$ for some $k\in {\bf R}$;

\vskip.03in
\item[(16)] a surface  defined by $(f,b,\psi,f): M\to \mathbb E^{m+1}_{s+1}$, where  $b$ is a unit time-like vector in $ \mathbb E^{m-5}_{s-1}$,  $\psi:M\to {\mathcal LC}^3_1\subset \mathbb E^4_1$ has constant mean curvature and $f$ is function $M$ satisfying $\Delta f=2k-2f$ for some $k\in {\bf R}$.
\end{enumerate}

\end{theorem} 

\section{Slant surfaces with parallel mean curvature vector in ${\bf C}^2_1$.}

  Slant surfaces with parallel mean curvature vector in the Lorentzian complex plane ${\bf C}^2_1$ were completely classified by the author together with Arslan, Carriazo and Mursthan  in \cite[2010]{Taiwan}.
 
 \begin{theorem}
  Let  $\Psi:M^2_1\to {\bf C}^2_1$ be a slant surface in the Lorentzian complex plane ${\bf C}^2_1$ with  parallel mean curvature vector. Then we have 
  \vskip.05in
  
 {\rm (A)} If $M^2_1$ is minimal in ${\bf C}^2_1$, then, up to rigid motions, $M^2_1$ is locally an open portion of one of the following three types of surfaces:

\vskip.05in
{\rm (A.1)} A totally geodesic slant plane;

\vskip.05in

{\rm (A.2)} A flat $\theta$-slant surface defined by
\begin{equation}\begin{aligned} \notag & \Psi(x,y)= \Big(x-\frac{\i y}{2}  \cosh\theta \coth \theta+(\i\sech\theta +\tanh \theta)K(y),
\\& \hskip.4in   x-y +\frac{\i y}{4}(\cosh 2 \theta-3)\csch\theta +(\i \sech\theta+\tanh\theta)K(y)
 \Big),\end{aligned}\end{equation}
where $K(y)$ is a non-constant function;

\vskip.05in

{\rm (A.3)} A non-flat $\theta$-slant surface defined by
\begin{equation}\begin{aligned} \notag&\hskip.1in  \Psi(x,y)=  \left( \int^y_0 \hskip-.04in  \frac{dy}{\sqrt{v'(y)}} - \frac{1+\i\sinh \theta}{2c} \int^x_0\hskip-.05in  \frac{u(x)dx}{\sqrt{u'(x)}}\right.
\\&\hskip.5in+\frac{ \i  \cosh \theta}{2b^{\frac{3}{2}}} \hskip-.03in \int^y_0 \hskip-.04in \frac{v(y)dy}{\sqrt{v'(y)}}+b^{\frac{3}{2}}c(\i \sech\theta-\tanh\theta)\hskip-.02in   \int^x_0\hskip-.08in  \frac{dx}{\sqrt{u'(x)}} , \\& \hskip.1in 
 \int^y_0\hskip-.05in \frac{dy}{\sqrt{v'(y)}}+\frac{1+\i \sinh\theta}{2c}\int^x_0  \hskip-.04in \frac{u(x)dx}{\sqrt{u'(x)}} 
\\&\hskip.5in\left.   -\frac{ \i \cosh \theta}{2b^{\frac{3}{2}}} \hskip-.03in \int^y_0 \hskip-.04in \frac{v(y)dy}{\sqrt{v'(y)}}+b^{\frac{3}{2}}c(\i \sech\theta-\tanh\theta)\hskip-.03in   \int^x_0\hskip-.08in  \frac{dx}{\sqrt{u'(x)}}\),
\end{aligned}\end{equation} 
where  $b,c$ are nonzero real numbers, and $u(x),v(y)$ are functions with $u'(x)>0,$ $v'(y)>0$ defined respectively on open intervals $I_1$ and $I_2$ containing $0$.
  
 \vskip.05in
 
{\rm (B)} If $M^2_1$ is non-minimal in ${\bf C}^2_1$,   then, up to rigid motions, $M^2_1$ is locally an open portion of one of the following nine types of flat slant surfaces in ${\bf C}^2_1:$
\vskip.05in

{\rm (B.1)}  A  Lagrangian surface defined by  $\Psi(x,y) =z(x) e^{\i a y}$,  where $a$ is a nonzero real number and $z$ is a null curve  lying in the light cone $\mathcal{LC}$  satisfying $\<\i z',z\>=a^{-1}$;
\vskip.05in

{\rm (B.2)}  A  Lagrangian surface defined by
\begin{equation}\begin{aligned}\notag &\Psi(x,y)=\Bigg(\frac{e^{\i cy}}{2c}\(2cx-\i+2\int^y_0 u(y)dy\)-\frac{1}{c}\int^y_0 e^{\i cy}u(y)dy,  
\\& \hskip.3in  \frac{e^{\i cy}}{2c}\(2cx+\i+2\int^y_0 u(y)dy\)-\frac{1}{c}\int^y_0 e^{\i cy}u(y)dy\Bigg),\end{aligned}\end{equation}
where $c$ is a nonzero real number and $u(y)$ is a nonzero real-valued function defined on an open interval  $I\ni 0$;

\vskip.05in
{\rm (B.3)}  A Lagrangian surface defined by
\begin{equation}\begin{aligned}  \notag &  \Psi(x,y)=\(\frac{x+ky}{\sqrt{2k}}, \frac{e^{2\i (x-ky)}}{2\sqrt{2k}}\), \end{aligned}\end{equation} where $k$ is a positive real number;

\vskip.05in
{\rm (B.4)}  A Lagrangian surface defined by
\begin{equation}\begin{aligned}  \notag &  \Psi(x,y)=\(\frac{e^{2\i (x+by)}}{2\sqrt{2b}},\frac{x-by}{\sqrt{2b}}\),\end{aligned}\end{equation} where $b$ is a positive real number;

\vskip.05in
{\rm (B.5)}  A Lagrangian surface defined by
\begin{equation}\begin{aligned}  \notag &  \Psi(x,y)=\frac{\sqrt{a}} {\sqrt{2b}} \(\frac{e^{\i (1+a^{-1})(ax+by)}} {a+1}, 
\frac{e^{\i (a^{-1}-1)(ax-by)}} {a-1}\),\end{aligned}\end{equation} 
where $a$ and $b$ are positive real numbers with $a\ne 1$;

\vskip.05in
{\rm (B.6)}  A Lagrangian surface defined by
\begin{equation}\begin{aligned}  \notag &  \Psi(x,y)=\frac{\sqrt{a}} {\sqrt{2k}} \(\frac{e^{\i (a^{-1}-1)(ax+ky)}} {a-1},
\frac{e^{\i (1+a^{-1})(ax-ky)}} {a+1}\),\end{aligned}\end{equation} 
where $a$ and $k$  are positive real numbers with $a\ne 1$;

\vskip.05in
{\rm (B.7)}  A Lagrangian surface defined by
\begin{equation}\begin{aligned}  \notag &  \Psi(x,y)={e^{(\i -a)x+(\i-a^{-1})by}}\(e^{2ax}-\frac{(a+\i)^4e^{2b^{-1}by}}{8b(1+a^2)^2}  , e^{2ax}-\frac{e^{2a^{-1}by}}{8b}\),\end{aligned}\end{equation} 
where $a$ is a positive real number and $b$ is a nonzero real number;

\vskip.05in
{\rm (B.8)}  A  proper slant surface with slant angle $\theta$ defined by
\begin{align}  \notag &\Psi(x,y)=z(x) \frac{(2y\sinh\theta-a\cosh\theta)^{\frac{1}{2}-\frac{\i}{2}\csch\theta}}{\sinh\theta -\i},\end{align} 
where $a$ is a real number  and $z(x)$ is a null curve  lying in the light cone $\mathcal{LC}$ which satisfies $\<z',\i z\>=\cosh^2\theta$;

\vskip.05in 
{\rm (B.9)}  A  proper slant surface with slant angle $\theta$  defined by
\begin{equation}\begin{aligned}  \notag & \hskip.3in  \Psi(x,y)=\Bigg( \sech^2\theta\int^y_0 u(y)(2y\sinh\theta-a\cosh\theta)^{\frac{3}{2}-\frac{\i}{2}\csch\theta}dy
  \\& +\hskip-.02in  (2y\sinh\theta\hskip-.02in -\hskip-.02in a\cosh\theta)^{\frac{1}{2}-\frac{\i}{2}\csch\theta}\times  \\& \hskip.4in  \(\hskip-.03in  x\hskip-.02in +\hskip-.02in \frac{\i}{2}\hskip-.02in 
  +\hskip-.02in  \sech^2\theta \hskip-.02in \int^y_0\hskip-.03in  (2y\sinh\theta-a\cosh\theta)u(y)dy\),
  \\&\hskip.3in \sech^2\theta\int^y_0 u(y)(2y\sinh\theta-a\cosh\theta)^{\frac{3}{2}-\frac{\i}{2}\csch\theta}dy
  \\& +\hskip-.02in  (2y\sinh\theta\hskip-.02in -\hskip-.02in a\cosh\theta)^{\frac{1}{2}-\frac{\i}{2}\csch\theta}\times  \\& \hskip.4in   \( \hskip-.03in x\hskip-.02in -\hskip-.02in \frac{\i}{2}\hskip-.02in 
  -\hskip-.02in  \sech^2\theta \hskip-.03in \int^y_0\hskip-.03in  (2y\sinh\theta-a\cosh\theta)u(y)dy\)\hskip-.05in \Bigg),\end{aligned}\end{equation}
where $a$ is a real number and $u(y)$ is a nonzero real-valued function defined on an open interval  $I\ni 0$.
 \end{theorem}

\section{Lorentz surfaces with parallel mean curvature vector in arbitrary pseudo-Euclidean space.}

Flat marginally trapped Lorentz surfaces with parallel mean curvature in  $\mathbb E^4_2$ were classified by the author in \cite[2008]{c3}. However, one of the cases in the classification given in \cite[2008]{c3} was missing, which was pointed out and  being added  in \cite[2010]{CY}. Based on this classification, the author and Garay classified  in \cite[2009]{CG}
all marginally trapped Lorentz surfaces with parallel mean curvature vector in  $\mathbb E^4_2$ (without the assumption on flatness on marginally trapped surfaces). 

Also based on the classification results of \cite[2008]{c3} and \cite[2010]{CY}, Hou and Yang classified in  \cite[2009]{HY} Lorentz surfaces \with in $\mathbb E^4_2$ (without condition on marginally trapped).

Finally, the author completely classified in \cite[2009]{Kyushu} all Lorentz surfaces \with in pseudo-Euclidean spaces with arbitrary codimension and index. More precisely, he proved  the following classification theorem.

\begin{theorem}\label{T:18.1} There are twenty-three families of Lorentz surfaces with parallel mean curvature in a pseudo-Euclidean $m$-space $\mathbb E^m_s$:

\vskip.04in
 {\rm (1)} A minimal surface of $\mathbb E^m_s$;

\vskip.04in
 {\rm (2)} A minimal surface of  a totally umbilical $S^{m-1}_s(c)\subset \mathbb E^m_s$ with $ c>0$;

\vskip.04in
 {\rm (3)} A minimal   surface of a totally umbilical $H^{m-1}_{s-1}(-c)\subset \mathbb E^m_s$ with $ c>0$;

\vskip.04in
{\rm (4)} A CMC surface of a totally geodesic $\mathbb E^3_1\subset \mathbb E^m_s$;

\vskip.04in
 {\rm (5)} A CMC surface of a totally geodesic $\mathbb E^3_2\subset \mathbb E^m_s$;

\vskip.04in
 {\rm (6)} A CMC surface of a totally umbilical $S^3_1(c)\subset\mathbb E^4_1\subset  \mathbb E^m_s$;

\vskip.04in
 {\rm (7)} A CMC surface of a totally umbilical $H^3_1(-c)\subset\mathbb E^4_2\subset  \mathbb E^m_s$;

\vskip.04in
 {\rm (8)} A CMC surface of a totally umbilical $S^3_2(c)\subset\mathbb E^4_1\subset  \mathbb E^m_s$;

\vskip.04in
 {\rm (9)} A CMC surface of a totally umbilical $H^3_2(-c)\subset\mathbb E^4_3\subset  \mathbb E^m_s$;

 {\rm (10)} A surface of  constant Gauss curvature lies in the light cone ${\mathcal LC}$ of a totally geodesic $\mathbb E^4_2\subset \mathbb E^m_s$;

\vskip.04in

  {\rm (11)} A surface lies in a totally geodesic $\mathbb E^5_2\subset \mathbb E^m_s$ defined by
$(f,\psi,f)$, where  $f$  is a function  satisfying $\Delta f=b\in {\bf R}$ and $\psi:M^2_1 \to \mathbb E^{3}_1$ is a CMC isometric immersion;
 \vskip.04in

  {\rm (12)} A surface lies in a totally geodesic $\mathbb E^5_3\subset \mathbb E^m_s$ defined by
$(f,\psi,f)$, where  $f$  is a function satisfying $\Delta f=b\in {\bf R}$ and $\psi:M^2_1 \to \mathbb E^{3}_2$ is a CMC isometric immersion;
 \vskip.04in
 
 {\rm (13)} A  surface of constant Gauss curvature $a^2>0$ defined by
$$L= \text{\Small$\frac{ w(y)}{x+y}$}-z(y),$$ where $w$ is a timelike curve  lying in the light cone ${\mathcal LC} \subset \mathbb E^m_s$  satisfying $$\<w,w''\>=-\<w',w'\>=\text{\Small$\frac{4}{a^2}$},\;\; \<w'',w''\>=- \text{\Small$\frac{2\d}{a}$},\; a>0$$ and $z$ is null curve satisfying $ z'=\frac{1}{2}(w''+a \d w)$ for some nonzero function $\d(y)$;

 \vskip.04in
{\rm (14)}  A   surface of constant Gauss curvature $-a^2<0$ defined by
$$L(x,y)=z(y)+ w(y)\tanh \Big(\!\text{\Small$\frac{a x+a y}{\sqrt{2}}$}\!\Big),$$ where $w$ is a timelike curve  lying in the light cone ${\mathcal LC} \subset \mathbb E^m_s$ satisfying $\<w',w'\>=-2$ and  
 $z$ is spacelike curve satisfying $$\<z',z'\>=2,\; \<z',w'\>=0,\; \<z',w\>=- \text{\Small$\frac{\sqrt{2}}{a}$},\;  z' =\text{\Small$\frac{1}{\sqrt{2}a}$}(a \d w-w''),\; \d(y)\ne 0; $$ 
 
\vskip.04in
{\rm (15)} A marginally trapped surface defined by $(f,\psi,f)$, where $\psi:M\to \mathbb E^{m-2}_{s-1}$ is a minimal immersion and $f$ is a function satisfying $\Delta f=c>0$,  $c\in {\bf R}$;

\vskip.04in
 {\rm (16)} A flat marginally trapped surface defined by $L=z(y) +x w(y)$, where $z,w$ are orthogonal null curves in $\mathbb E^m_s$ satisfying $\<w(y),z'(y)\>=-1$ and $w''(y)=v(y)w(y)$ for some nonzero function $v$, and  $w$ lies in the light cone;

\vskip.04in
 {\rm (17)}  A  flat marginally trapped surface surface defined by $L=z(x)y+w(x),$
where $z,w$ are null curves in a totally geodesic $ \mathbb E^4_2\subset \mathbb E^m_s$ and
$z$ lies in the light-cone $\mathcal{LC}$ which  satisfies $\<z',w'\>=0,\<z,w'\>=-1$ and $z''=-\beta z$ for a function $\beta(x)$.

\vskip.04in
{\rm (18)}  A flat marginally trapped surface  defined by  $L=z(y)\cos a x+ w(y)\sin a x,$ where  $z,w$ are null curves lying in  the light cone  $\mathcal{LC}$ of a  totally geodesic $ \mathbb E^4_2\subset \mathbb E^m_s$ such that $\<z,w\>=z''+\delta z=w''+\delta w=0$ and $\<z,w'\>=a^{-1}$ for some non-constant function $\delta(y)$ and positive number $a$.

\vskip.04in
 {\rm (19)} A flat marginally trapped surface  defined by $L=z(y)\cosh ax+w(y)\sinh ax,$
 where  $z,w$ are null curves lying in  the light cone  $\mathcal{LC}$ of a  totally geodesic $ \mathbb E^4_2\subset \mathbb E^m_s$ such that $\<z,w\>=z''+\delta z=w''+\delta w=0$ and $\<z,w'\>=a^{-1}$ for some non-constant function $\delta(y)$ and positive number $a$.

\vskip.04in
 {\rm (20)} A flat marginally trapped surface lies in a totally geodesic $\mathbb E^4_2\subset \mathbb E^m_s$ defined by
$$\text{\Small$ \frac{1}{\sqrt{2}}$}\big(x y+f(x)+k(y),x+y,x-y, xy+f(x)+k(y)\big)$$ for some functions $f,k$.

\vskip.04in
{\rm (21)}  A flat marginally trapped surface lies in a totally geodesic $\mathbb E^4_2\subset \mathbb E^m_s$ defined by
\begin{equation}\begin{aligned}\notag&\text{\Small$\frac{1}{2ab}$}\Big( 2ab\cos ax\cos by-\sin ax\sin by, 2ab\cos ax\sin by+\sin ax\cos by,\\& \hskip-.1in 2ab\cos ax\cos by+\sin ax\sin by,2ab\cos ax\sin by-\sin ax\cos by\Big),\, a,b>0;\end{aligned}\end{equation}

\vskip.04in
{\rm (22)}  A flat marginally trapped surface lies in a totally geodesic $\mathbb E^4_2\subset \mathbb E^m_s$ defined by
\begin{equation}\begin{aligned}\notag&\text{\Small$\frac{1}{2ab}$}\Big( 2ab\cos ax\cosh by-\sin ax\sinh by, 2ab\cos ax\sinh by+\sin ax\cosh by,\\& \hskip-.1in 2ab\cos ax\cosh by+\sin ax\sinh by,\\&\hskip.4in 2ab\cos ax\sinh by-\sin ax\cosh by\Big),\, a,b>0;\end{aligned}\end{equation} 

\vskip.04in
 {\rm (23)} A flat marginally trapped surface lies in a totally geodesic $\mathbb E^4_2\subset \mathbb E^m_s$ defined by
\begin{equation}\begin{aligned}\notag&\text{\Small$\frac{1}{2ab}$}\Big( 2ab\cosh ax\cosh by-\sinh ax\sinh by, 2ab\cosh ax\sinh by\\&\hskip.2in  +\sinh ax\cosh by, 2ab\cosh ax\cosh by+ \sinh ax\sinh by, \\&\hskip.4in 2ab\cosh ax\sinh by-\sinh ax\cosh by\Big),\, a,b>0.\end{aligned}\end{equation} 
\vskip.04in

Conversely, each  Lorentz surface with parallel mean curvature in  $\mathbb E^m_s$  is congruent to a surface given by one of the twenty-three families described above.
\end{theorem}

\begin{remark} Let  $\psi:M^2_1\to S^m_s(1)$ (resp. $\psi:M^2_1\to H^m_s(-1))$ be an isometric immersion of a Lorentz surface $M^2_1$. Then the composition $\iota\circ \psi:M\to S^m_s(1)\subset \mathbb E^{m+1}_s$ (resp. $\iota\circ \psi : M^2_1\to H^{m}_s(-1)\subset \mathbb E^{m+1}_{s+1}$) via (1.2) (resp. via (1.3)) has parallel mean curvature vector if and only if $\psi$ has parallel mean curvature vector. Consequently, to classify Lorentz surfaces in $S^m_s(1)$ with parallel mean curvature vector, we only need to determine those Lorentz surfaces $M^2_1$ in $\mathbb E^{m+1}_s$ with parallel mean curvature vector such that $M^2_1$ lies in a $S^m_s(1)\subset \mathbb E^{m+1}_s$ via (1.2). Similar method applies to Lorentz surfaces in $H^m_s(-1)$ with parallel mean curvature vector.
\end{remark}

\begin{remark} Space-like surfaces with parallel second fundamental form in indefinite space forms have parallel mean curvature vector. Such surfaces have been completely classified by the author in \cite[2010]{JGP}.
\end{remark}

\begin{remark}Lorentz surfaces with parallel second fundamental form in indefinite space forms have parallel mean curvature vector. Such surfaces have been completely classified in a series of recent articles \cite{1,2,3,4,5,cv}.

\end{remark}

  \end{document}